\documentclass[12pt,draftcls,onecolumn]{IEEEtran}

\ifCLASSINFOpdf
\else
\fi
\usepackage{graphicx}
\usepackage{epstopdf}
\usepackage{epsfig}
\usepackage{amsfonts}
\usepackage{amstext}
\usepackage{color}
\usepackage{empheq}
\usepackage{graphicx,colordvi}
\usepackage{amssymb,amsmath,amsopn,amsfonts,graphicx}
\usepackage{CJK}
\usepackage{float}

\newtheorem{lemma}{\bf Lemma}[section]
\newtheorem{theorem}{\bf Theorem}[section]
\newtheorem{definition}{\bf Definition}[section]

\newtheorem{remark}{\bf Remark}[section]

\newtheorem{proposition}{\bf Proposition}[section]

\usepackage[justification=centering]{caption}

\begin{document}
%\title{Bare Demo of IEEEtran.cls\\ for IEEE Journals}

\title{Indefinite Linear Quadratic Mean Field Social Control Problems with Multiplicative Noise}

\author{
	Bingchang Wang, \emph{Member, IEEE,}
	\thanks{This work was supported by the National Natural Science Foundation of China under Grants 61773241, 61573221 and
61633014.}%, and in part by Taishan Scholar Construction Engineering by Shandong Government.}
	\thanks{Bingchang Wang is with the School of Control Science and Engineering, Shandong University, Jinan 250061, P. R. China. (e-mail: bcwang@sdu.edu.cn) }
and Huanshui Zhang, \emph{Senior Member, IEEE}

	\thanks{Huanshui Zhang is with the School of Control Science and Engineering, Shandong University, Jinan 250061, P. R. China. (e-mail: hszhang@sdu.edu.cn) }

}


\markboth{Journal of \LaTeX\ Class Files}%
{Shell \MakeLowercase{\textit{et al.}}: Bare Demo of IEEEtran.cls for IEEE Journals}



\maketitle

\begin{abstract}
This paper studies uniform stabilization and social optimality for linear quadratic (LQ) mean field control problems with multiplicative noise,
where agents are coupled via dynamics and individual costs. The state and control weights in cost functionals are not limited to be positive semi-definite.
This leads to an \emph{indefinite} LQ mean field control problem, which may still be well-posed due to deep nature of multiplicative noise.  %(teams)
We first obtain a set of forward-backward stochastic differential equations (FBSDEs) from variational analysis, and construct a feedback control by decoupling the FBSDEs. By using solutions to two Riccati equations, we design a set of decentralized control laws, which is further shown to be asymptotically social optimal. Some equivalent conditions are given for uniform stabilization of the systems with the help of linear matrix inequalities. A numerical example is given to illustrate the effectiveness of the proposed control laws.

%For the game problem, we first design a set of decentralized control from variational analysis, and then show that such set of decentralized control constitute an asymptotic Nash equilibrium by exploiting the stabilizing solution of a nonsymmetric Riccati equation.

%It is verified that the proposed decentralized control laws are equivalent to the feedback strategies of mean field control in previous works. This may illustrate the relationship between open-loop and feedback solutions of mean field control (games).
%Some numerical examples are given to illustrate the effectiveness of the proposed control laws.
%Finally, we provide
\end{abstract}
%
%
%% Note that keywords are not normally used for peerreview papers.
%
%
%
%
%
%
%% For peer review papers, you can put extra information on the cover
%% page as needed:
%% \ifCLASSOPTIONpeerreview
%% \begin{center} \bfseries EDICS Category: 3-BBND \end{center}
%% \fi
%%
%% For peerreview papers, this IEEEtran command inserts a page break and
%% creates the second title. It will be ignored for other modes.
%\IEEEpeerreviewmaketitle
%
\section{Introduction}

\subsection{Background and motivation}
The topic of mean field games and control
has drawn increasing attention in many disciplines including system control, applied mathematics and economics \cite{BFY13}, \cite{C14}, \cite{GS13}. A mean field game involves a very large number of small interacting players. While the influence of each player is negligible, the impact of the overall population is significant. By combining mean field approximations and individual best response,
the dimensionality difficulty can be overcome.  Mean field games and control
have found wide applications, including smart grids \cite{MCH13}, \cite{CBM15}, finance, economics \cite{GLL11}, \cite{CS14}, \cite{WH15}, \cite{HN16a}, %operations research %\cite{LW11}, \cite{AJW13},
%%\cite{LM14}, %\cite{SMJ15},
and social networks %\cite{LT15},
\cite{BTN16}, \cite{LW19}, etc.
%By identifying a consistency relationship between the individual's best response  and the mass (population macroscopic) behavior, one may
% obtain a fixed-point equation to specify the mean field.  This procedure leads to a set of  decentralized  strategies as an $\epsilon$-Nash equilibrium for the actual model with a large but finite population.
%Mean field models have widely appeared in economics, engineering, and social sciences, for instance, output planning in large markets \cite{L84}, dynamic advertising competition \cite{E09}, wireless communication networks \cite{HCM03}, and voluntary vaccination games \cite{B04}.

Depending on the state-cost setup of a mean field game, it can be classified into linear-quadratic (LQ) type and more general nonlinear type. By now, the LQ type has been commonly adopted in mean field studies because of its analytical tractability and close connection to practical applications. In this aspect, some relevant works include \cite{HCM07}, \cite{LZ08}, \cite{WZ13}, \cite{BSYY16}, \cite{MB17}.  %HWW16,
% there are also a large body of works on nonlinear models \cite{HMC06, YMMS12, CD13, DH15}.
Huang \emph{et al.} developed the Nash certainty equivalence (NCE) based on the fixed-point method and designed an $\epsilon$-Nash equilibrium for LQ games with discount costs \cite{HCM07}. The NCE approach was then applied to the (general) cases with stochastic ergodic costs \cite{LZ08} and with Markov jump parameters \cite{WZ13}, respectively. %HCM03
%Carmona and Delarue considered mean field
%games via the stochastic maximum principle, %They assumed that the dynamics are affine in states and controls, and the costs are convex,
%and obtained the $\epsilon$-Nash equilibrium by Schauder's Theorem
 The works \cite{CD13}, \cite{BSYY16} employed the adjoint equation approach and the fixed-point
theorem to obtain sufficient conditions for the existence of the
equilibrium strategy over a finite horizon.
%Lasry and Lions independently introduced the general model of mean field games and studied the well-posed problem of the limiting partial differential equations \cite{LL07}.
For other aspects of mean field games, readers are referred to \cite{HMC06}, \cite{LL07}, \cite{CD13} for nonlinear mean field games, \cite{weintraub2008markov} for oblivious equilibrium in dynamic games, \cite{H10}, \cite{WZ12}, \cite{WZ14} for mean field games with major players, \cite{HH16}, \cite{MB17} for robust mean field games.
%\cite{YMMS12} for the game of the synchronization in nonlinear oscillators, \cite{BSY13, DH15} for time-inconsistent mean field games, \cite{weintraub2008markov} for the notion of oblivious equilibrium in dynamic games. %\cite{CD13} for probabilistic analysis of mean field games.
%\cite{BSY13, DH15} for time-inconsistent mean field games.

Apart from noncooperative games, team optimization forms another research branch for studying cooperative behavior among multiple decision makers. In particular, social optima in mean field models with weak coupling have drawn more research interests. By social optimization, all players in a large population system (endowed with some weak-coupling structure in either cost or dynamics) will cooperate to minimize a common social cost--the sum of individual costs. Accordingly, we formulate a type of team decision problem \cite{R62}. Different from Nash games, all the agents in a team problem
are cooperative and share the same cost functional, although
they may have different information sets \cite{H80}.
Huang \emph{et al.} considered social optima in mean field LQ control, and provided an asymptotic team-optimal solution \cite{HCM12}. Wang and Zhang \cite{WZ17} investigated a mean field social optimal problem where the Markov jump parameter appears as a common source of randomness. %The work \cite{HN16} designed socially optimal strategies for mixed games by analyzing forward-backward stochastic differential equations (FBSDEs).
For further literature on social control, for instance, see \cite{HN16} for social optima in mixed games, \cite{AM15} for team-optimal control with finite population and partial information,
and \cite{SNM18} for the dynamic %cooperative
collective choice by finding a social optimum.
%social optima in static mean field games.

Concerned with mean field games and control for stochastic systems, most existing literature focused on the case with additive noise (i.e., the intensity of noise is independent of the state). Sometimes, such kind of noise is not sufficient to depict practical situations.
Alternatively, multiplicative noise is another realistic description for stochastic disturbance. Mean field control with multiplicative noise has attracted much attention due to its wide applications in engineering,
economics, and etc \cite{CD18}, \cite{HN16a}, \cite{WNZ19}, \cite{WZ19}.
This paper investigates %social optima for linear quadratic mean field control systems,
uniform stabilization and social optimality for mean field LQ control systems with multiplicative noises,
where subsystems are coupled via both dynamics and individual costs. The intensities of multiplicative noises depends on both system states and control inputs. %(teams)
The state weight $Q$ and control weight $R$ in the cost functional are not limited to be positive semi-definite. %In our model,
 In fact, an \emph{indefinite} $Q$ or $R$ may naturally occur in a wide class of practical problems, including production adjustment \cite{WH15},  uncertain systems \cite{HH16}, and portfolio selection \cite{ZL00}.
This problem leads to generalized %stochastic
Riccati equations, which is essentially different from the classical Riccati equation %since the problem is indefinite and involves multiplicative noise.
due to indefinite weights and multiplicative noise appearing in the problem.

\subsection{Challenge and main contributions}
{Most previous results on mean field games and control were given by virtue of the fixed-point analysis \cite{HCM07}, \cite{LZ08}, \cite{HCM12}, \cite{CD13}, \cite{BSYY16}, \cite{WZ17}.
	However, the fixed-point assumption may be
not easy to tackle, particularly for high-dimensional systems.
In this paper, we %break away from the fixed-point method and
solve the problem by decoupling
	forward-backward stochastic differential equations (FBSDEs) instead of fixed-point analysis. In recent years, some substantial progress for the optimal LQ control has been made by solving the FBSDEs. See \cite{Y13}, \cite{ZX17}, \cite{ZQ16}, \cite{NLZ16}, \cite{W18}, \cite{QZW19} for details.}

%We first consider
For the finite-horizon mean field LQ control problem, we first obtain a set of FBSDEs by examining the social cost variation, and give a centralized feedback control by decoupling the FBSDEs. Applying mean field approximations, we design decentralized control laws. By exploiting the uniform convexity property of the optimal control problem, we further show that the decentralized controls have asymptotic social optimality.
For the infinite-horizon case, we construct a set of decentralized control laws by using solutions of two Riccati equations, and further show decentralized controls are asymptotically social optimal. Some equivalent conditions are further given for uniform stabilization of all the subsystems with the help of linear matrix inequalities. %and Hamiltonian matrix eigenvalues. %in different cases.when the state weight $Q$ is semi-positive definite or only symmetric.
% 乘性噪声 一致镇定和平均场逼近的证明更加困难

For the mean field control systems with multiplicative noise, it is more difficult to show the uniform stabilization of all the subsystems than the case with additive noise. %and the consistency of mean field approximations.
Due to the appearance of multiplicative noise, the approximation error between population state average $\hat{x}^{(N)}$ and aggregate effect $\bar{x}$ relies on the states of all the agents while the mean square of the state $\hat{x}_i$ conversely depends on the approximation error. Thus, we need to analyze jointly the approximation error and states of all the agents. By tackling the corresponding integral inequalities, we obtain that all the subsystems are uniformly stabilizable and the mean field approximation is consistent.
%we need to exploit further the relation between population state average and aggregate effect. joint estimation of state and approximation error
Moreover, %indefinite nature of the problem
 since the weights $Q$ and $R$ in the cost functional are \emph{indefinite}, the prior boundedness of the state is not implied directly by the finiteness of the cost, which brings about extra difficulty to show the social optimality of decentralized control. %Thus, it is more difficult to show the social optimality of decentralized control.
 Here we first obtain the prior upper bounds of states and controls by exploiting the uniform convexity property of the problem, and further prove that decentralized strategies have asymptotic social optimality by perturbation analysis.
%Finally, some numerical examples are given to illustrate the effectiveness of the proposed control laws.

The main contributions of the paper are summarized as follows.
\begin{itemize}
 \item For the finite-horizon problem, we first obtain necessary and sufficient existence conditions of centralized optimal control based on FBSDEs, and then design a feedback-type decentralized control by decoupling FBSDEs and applying mean field approximations.

%(ii)
\item By exploiting the uniform convexity of the problem, the decentralized control laws are shown to have asymptotic social optimality.

%(iii)
\item The necessary and sufficient conditions are given for uniform stabilization of the systems by virtue of the system's observability and linear matrix inequalities.

    \item An explicit expression of the asymptotic average social cost is given in terms of the
solutions of two Riccati equations.
\end{itemize}

%(iv) For the game problem, we show that the decentralized control laws constitute an $\varepsilon$-Nash equilibrium by exploiting the stabilizing solution of a nonsymmetric Riccati equation.
%
%(v) It is under nonconservative assumptions that we obtain the asymptotically optimal decentralized control, and such control laws are shown to be equivalent to the feedback strategies given by the fixed-point method in previous works \cite{HCM07, HCM12}.
\subsection{Organization and notation}

The organization of the paper is as follows. In Section II, the indefinite LQ mean field social control problem is formulated.  In Section III, we first construct a set of decentralized control laws for the finite-horizon case, and then show its asymptotic social optimality. In Section IV, we design asymptotically optimal control for the infinite-horizon case and further give some equivalent conditions of uniform stabilization. In Section V, we give the value of asymptotic average optimal social cost.
In Section VI, a numerical example is provided to show the effectiveness of the proposed controls.
Section VII concludes the paper.

The following notation will be used throughout this paper.
Denote by $\|\cdot\|$ the Euclidean vector norm or matrix spectral norm, and $\otimes$ the Kronecker product. For a vector $z$ and a matrix $Q$, $\|z\|_Q^2= z^TQz$; $Q>0$ ($Q\geq0$) means that the matrix $Q$ is positive definite (positive semi-definite).
$Q^{\dag}$ is the Moore-Penrose pseudoinverse\footnote{$Q^{\dag}$ is a unique matrix satisfying
$QQ^{\dag}Q=Q^{\dag}, Q^{\dag}QQ^{\dag}=Q, (Q^{\dag}Q)^T=Q^{\dag}Q$, and $(QQ^{\dag})^T=QQ^{\dag}.$} of the matrix $Q$, $\mathcal{R}(Q)$ denotes the range of a matrix (or an operator) $Q$, and $\hbox{ker}(Q)$ is the kernel of $Q$.
%For a matrix function $R(t)$, $R(t)\gg 0$ means $R(t)\geq\lambda I$ for some $\lambda>0$ and any $t$.
For two vectors $x,y$, $\langle x,y\rangle=x^Ty$.
%$C([0,\infty),\mathbb{R}^ k)$ is the space of all $\mathbb{R}^ k$-valued continuous functions defined on $[0,\infty)$, and
$L^2([0,\infty), \mathbb{R}^k)$ is given by $\{f: [0,\infty)\to \mathbb{R}^k|\int_0^{\infty}\|f(t)\|^2dt<\infty  \}.$
$L^2_{\mathcal F}(0, T; \mathbb{R}^k)$ is the space of all $\mathcal{F}_t$-adapted $\mathbb{R}^k$-valued processes $x(\cdot)$ such that
$\mathbb{E}\int_0^T\|x(t)\|^2dt<\infty$.
%
%For two sequences $\{a_n, n = 0, 1, \cdots\}$ and $\{b_n, n = 0, 1,  \cdots\}$, $a_n = O(b_n)$ denotes that
%$\limsup_{n\to\infty}|{a_n}/{b_n}|\leq c$, and $a_n = o(b_n)$ denotes
%$\limsup_{n\to\infty}|{a_n}/{b_n}|=0$.
For convenience of  presentation, we use $c, c_1,c_2,\cdots$ to
denote generic positive constants, which may vary from place to place.

\section{Problem Description}%{Mean Field LQ Social Control}
\label{sec2.3.1}

Consider a large population systems with $N$ agents. Agent $i$ evolves by the following stochastic differential equation:
\begin{align}\label{eq1}
  dx_i(t) = &[Ax_i(t)+Bu_i(t)+Gx^{(N)}(t)+f(t)]dt\cr
  &+[Cx_i(t)+Du_i(t)+\sigma(t)] dW_i(t),\ 1\leq i\leq N,
\end{align}
where $x_i\in
\mathbb{R}^n$ and $u_i\in\mathbb{R}^r$ are the state and input of the $i$th agent.  $A,B,G, C,D$ are constant matrices with appropriate dimensions. $x^{(N)}(t)=\frac{1}{N}\sum_{j=1}^Nx_j(t)$, $f, \sigma\in L^2([0,\infty), \mathbb{R}^n)$.
$\{W_i(t),1\leq i\leq N\}$ are a sequence of independent $1$-dimensional Brownian motions on a complete
filtered probability space $(\Omega,
\mathcal F, \{\mathcal F_t\}_{0\leq t\leq T}, \mathbb{P})$.
The cost functional of agent $i$ is given by
\begin{equation}\label{eq2}
J_i(u)= \mathbb{E}\int_0^{\infty}\Big\{\big\|x_i(t)
-\Gamma x^{(N)}(t)-\eta(t)\big\|^2_{Q}
+\|u_i(t)\|^2_{R}\Big\}dt,
\end{equation}
where $Q,R, \Gamma\in \mathbb{R}^{n\times n}$ are constant matrices, and $\eta\in L^2([0,\infty), \mathbb{R}^n)$. $Q$ and $R$ are symmetric (generally indefinite). % and $R>0$.
Denote $u=\{u_1,
\ldots,u_{i}, \ldots, u_N\}$.
The decentralized control set is given by
\begin{equation*}
{\mathcal U}_{d} =\Big\{(u_1,\cdots,u_N)\ \big|\
u_i(t)\ \hbox{is adapted to}
\ \sigma(x_i(s),0\leq s\leq t),
\mathbb{E}\int_0^{\infty}\|x_i(t)\|^2dt<\infty,\forall i\Big\}.
%{\cal U}_{d,i} =\Big\{u_i\ \big|\ &
%u_i(t)\ \hbox{is adapted to}\  \sigma(x_i(s),0\leq s\leq t),  \mathbb{E}\int_0^{\infty}\|x_i(t)\|^2dt<\infty\Big\}.
\end{equation*}
For comparison, define the centralized control set as
\begin{equation*}
  {\mathcal U}_{c} =\Big\{(u_1,\cdots,u_N) \big|\
u_i(t)\ \hbox{is adapted to} \ {\mathcal F}_t,
\mathbb{E}\int_0^{\infty}\|x_i(t)\|^2dt<\infty, \forall i\Big\},
\end{equation*}
%and $
 % {\cal U}_{c} =\big\{(u_1,\cdots,u_N) \big|\
%u_i\ \hbox{is adapted to}\ {\cal U}_{c,i}$\big\},
where ${\mathcal F}_t\stackrel{\Delta}{=} \sigma\{\bigcup_{i=1}^N{\mathcal F}_t^{i}\}$ and ${\mathcal F}_t^{i} = \sigma(x_{i}(0), W_i(s),0\leq s\leq t), i=1,\cdots,N$.

In this paper, we mainly study the following problem.

\textbf{(P0)} Seek a set of decentralized control laws to optimize social cost
for the system (\ref{eq1})-(\ref{eq2}), i.e.,
$\inf_{u\in {\mathcal U}_{d}}J_{\rm soc}(u),$
where $$J_{\rm soc}(u)=\sum_{i=1}^NJ_i(u).$$
%$$\begin{aligned}
%{\cal U}_{c,i} &=\Big\{u_i|\
%  u_i(t)\in \sigma(\bigcup_{i=1}^n{\mathcal F}_t^{i})\Big\}.
%\end{aligned}
%$$
%Problem (P1): maximize ${J}_{soc}$ over $f \in {\mathcal U}_{c}$.

We first make the assumption on the initial values of agents' states.

\textbf{A1)} $x_i(0), i=1,...,N$ are mutually independent and have the same mathematical expectation. $x_{i}(0)=x_{i0}$, $\mathbb{E}x_i(0)=\bar{x}_0$, $i=1,\cdots,N$. There exists a constant $c_0$ (independent of $N$) such that $\max_{1\leq i \leq N}\mathbb{E}\|x_i(0)\|^2<c_0$. Furthermore, $\{x_i(0), i=1,...,N\}$ and
$\{W_i, i=1,...,N\}$ are independent of each other.

\begin{remark}
Since the weights $Q$ and $R$ are indefinite, Problem (P0) is called an indefinite LQ mean field social control problem.
Due to the indefiniteness of $Q$ and $R$, the convexity may be lost, and the problem may have no solutions.
Thus, we need to discuss the convexity of Problem (P0), which is related to the generalized Riccati equation.
\end{remark}

To facilitate the discussion for the convexity of Problem (P0), we write the problem in a high-dimensional form.

Let  $\textbf{x}=(x_1^T,\cdots,x_N^T)^T$, $\textbf{u}=(u_1^T,\cdots,u_N^T)^T$, $\textbf{1}=(1,\cdots,1)^T$, ${\mathbf{\sigma}_i}=(0,\cdots,0,\sigma^T, 0,\cdots,0)^T$,
$\check{\textbf{A}}=diag(A, \cdots,A)+\frac{1}{N}(\textbf{11}^T\otimes G)$, $\textbf{B}=diag(B, \cdots,B)$, $\textbf{C}_i=diag(0, \cdots,0,C,0,\cdots,0)$, $\textbf{D}_i=diag(0, \cdots,0,D,0,\cdots,0)$, and $\textbf{R}=diag(R, \cdots,R)$.
With the above notations, the dynamics of all agents can be written in the more compact form:
\begin{equation*}
d\textbf{x}(t) = \big[\check{\textbf{A}}\textbf{x}(t)+\textbf{Bu}(t)+ \textbf{1}\otimes f(t)\big]dt
+\sum_{i=1}^N[\textbf{C}_i\textbf{x}(t)+\textbf{D}_i\textbf{u}(t)+{\bf{\sigma}}_i(t)] dW_i(t).
\end{equation*}
Also, denote
$$
\left\{ \begin{array}{l}
Q_{\Gamma}\stackrel{\Delta}{=}\Gamma^TQ+Q\Gamma-\Gamma^TQ\Gamma\\
\bar{\eta}\stackrel{\Delta}{=}Q\eta-\Gamma^T Q\eta
\end{array}
\right.$$
  By rearranging the integrand of ${J}_{soc}$, we have
\begin{equation}\label{eq4}
 {J}_{soc}%=&\frac{1}{2}\sum_{i=1}^N\mathbb{E}\int_0^{T}\Big(-\textbf{x}^T{\bf{\Gamma}}_i^T{Q}{\bf{\Gamma}}_i\textbf{x}+2\eta^TQ{\bf{\Gamma}}_i\textbf{x}
 %+\sigma_0^TR_0\sigma_0\Big)dt\cr
 =\mathbb{E}\int_0^{\infty}\Big(\|\textbf{x}(t)\|_{\bar{\textbf{Q}}}^2-2(\textbf{1}\otimes \bar{\eta}(t))^T\textbf{x}(t)
 +N\|\eta(t)\|^2+\|\textbf{u}(t)\|_{\textbf{R}}^2\Big)dt,
\end{equation}
where %$\bar{\bf{\eta}}=\textbf{1}\otimes \bar{\eta}$, and
 $\bar{\textbf{Q}}=(\bar{Q}_{ij})$ is given by
$$
\bar{Q}_{ii}=Q-Q_{\Gamma}/N,\ \bar{Q}_{ij}=-Q_{\Gamma}/N,\ 1\leq i\neq j\leq N.
$$

 \begin{remark}\label{rem1}
{Hereafter, we may exchange the usage of notation $u=(u_1, \cdots, u_{N}) \in \mathbb{R}^{r \times N}$ and $\textbf{u}=(u_1^T,\cdots,u_N^T)^{T} \in \mathbb{R}^{rN}$. Both notations represent the control laws among all agents, but only differ in their formations.}
 \end{remark}

\section{Mean Field LQ Social Control over a Finite Horizon}\label{sec3}
%Denote  {\cal U}_{c}
For the convenience of design, we first consider the following finite-horizon problem.\\
$$\textbf{(P1)}\inf_{u\in L^2_{{\mathcal F}}(0, T; \mathbb{R}^{nr})} J_{\rm soc}^{\rm F}(u),%=\inf_{u_i\in {\cal U}_{d,i}}\sum_{i=1}^NJ_{i}^{\rm F}(u),
$$
where $J_{\rm soc}^{\rm F}(u)=\sum_{i=1}^NJ_{i}^{\rm F}(u)$ and
\begin{equation}\label{eq3}
J_{i}^{\rm F}(u)=\mathbb{E}\int_0^{T}\Big\{\big\|x_i(t)
-\Gamma x^{(N)}(t)-\eta(t)\big\|^2_{Q}+\|u_i(t)\|^2_{R}\Big\}dt
+\mathbb{E}\|x_i(T)-\Gamma_0x^{(N)}(T)-\eta_0\|_{H}^2.
\end{equation}

We now provide some equivalent conditions for the convexity of Problem (P1). Denote
$$
H_{\Gamma_0}\stackrel{\Delta}{=}\Gamma^T_0H+H\Gamma_0-\Gamma^T_0H\Gamma_0,\
\bar{\eta}_0\stackrel{\Delta}{=}H\eta_0-\Gamma^T_0 H\eta_0.$$

\begin{proposition}\label{prop1}
 Problem (P1) is convex in $u$ if and only if
for any $u_i\in L^2_{{\mathcal F}}(0, T; \mathbb{R}^r)$, $i=1,\cdots,N$,
$$
  \sum_{i=1}^N\mathbb{E}\!\int_0^T\!\Big\{\big\|y_i-\Gamma y^{(N)}\big\|^2_{Q}+\|u_i\|^2_{R}\Big\}dt
  +\sum_{i=1}^N\mathbb{E}\|y_i(T)-\Gamma_0y^{(N)}(T)\|_{H}^2\geq 0,
$$
where $y^{(N)}=\sum_{j=1}^Ny_j/N$ and $y_i$ satisfies
\begin{align}\label{eq4aa}
   dy_i=(A y_i+Gy^{(N)}+Bu_i)dt+(Cy_i+Du_i)dW_i,\ y_i(0)=0,\quad i=1,\cdots,N.
\end{align}
\end{proposition}

\emph{Proof.} The proof is similar to \cite{HH16}, \cite{LZ99}. $\hfill \Box$
%
%\emph{Proof.} Let $x_i$ and $\acute{x}_i$ be the state processes of agent $i$ %(\ref{eq1})
%with the control $v$ and $\acute{v}$,
%respectively.
% Take any $\lambda_1\in [0, 1]$ and let $\lambda_2= 1-\lambda_1$.
%Then
%\begin{align*}
%&\lambda_1  J_{\rm soc}^{\rm F} (v) +\lambda_2 J_{\rm soc}^{\rm F}  (\acute{v}) - J_{\rm soc}^{\rm F}(\lambda_1 v+\lambda_2 \acute{v})  \\
% =& \lambda_1\lambda_2 \mathbb{E}\int_0^T\left\{ \|x_i-\acute{x}_i -\Gamma(x^{N}-\acute{x}^{N})\|_Q^2 + \|u_i-\acute{u}_i\|^2_R\right\} dt.
%\end{align*}
%Denote $u= v-\acute{v}$, and $y_i=x_i-\acute{x}_i$. Thus, $y_i$ satisfies
%(\ref{eq4aa}). %Hence
%%$$
%%\lambda_1 \bar J'_{i} (u_i,  f) +\lambda_2 \bar J'_{i} (u_i,  f') - \bar J'_{i} (u_i,  \lambda_1 f+\lambda_2 f')
%% =\lambda_1\lambda_2\bar J_{i}''(g)
%%$$
%By the definition of the convexity, the lemma follows.  $\hfill \Box$

\begin{proposition}
  \label{prop2}
  The following statements are equivalent:

(i) Problem (P1) is uniformly convex in $u$;

(ii) For any $u_i\in L^2_{{\mathcal F}}(0, T; \mathbb{R}^r)$, $i=1,\cdots,N$, there exists a constant $\gamma>0$ such that
$$  \begin{aligned}
  &\sum_{i=1}^N\mathbb{E}\int_0^T\Big\{\big\|y_i-\Gamma y^{(N)}\big\|^2_{Q}+\|u_i\|^2_{R}\Big\}dt
  +\sum_{i=1}^N\mathbb{E}\|y_i(T)-\Gamma_0y^{(N)}(T)\|_{H}^2\cr
  \geq& \gamma\sum_{i=1}^N \mathbb{E}\int_0^T\|u_i\|^2dt,
\end{aligned} $$
(iii) The equation
\begin{align}\label{eq5}
  \dot{\textbf{P}}+
\check{\textbf{A}}^T\textbf{P}+\textbf{P}\check{\textbf{A}}+\sum_{i=1}^N\textbf{C}_i^T\textbf{P}\textbf{C}_i+\bar{\textbf{Q}}-
\Big(\textbf{B}^T\textbf{P}+\sum_{i=1}^N\textbf{D}_i^T\textbf{P}\textbf{C}_i\Big)^T
\mathbf{\Upsilon}^{\dag}
\Big(\textbf{B}^T\textbf{P}+\sum_{i=1}^N\textbf{D}_i^T\textbf{P}\textbf{C}_i\Big)=0,
\end{align}
with $\textbf{P}(T)=\bar{\textbf{H}}$ admits a solution such that $\mathbf{\Upsilon}=\textbf{R}
+\sum_{i=1}^N\textbf{D}_i^T\textbf{P}\textbf{D}_i\geq 0$ and $\mathcal{R}\Big(\textbf{B}^T\textbf{P}+\sum_{i=1}^N\textbf{D}_i^T\textbf{P}\textbf{C}_i\Big)
\subseteq\mathcal{R}(\mathbf{\Upsilon})$,
where %$\bar{\bf{\eta}}=\textbf{1}\otimes \bar{\eta}$, and
 $\bar{\textbf{H}}=(\bar{H}_{ij})$ is given by
$$
\bar{H}_{ii}=H-H_{\Gamma_0}/N,\ \bar{H}_{ij}=-H_{\Gamma_0}/N,\ 1\leq i\neq j\leq N.
$$
\end{proposition}

\emph{Proof.}
(i)$\Leftrightarrow$(ii) is implied from \cite{HH16}, \cite{LZ99}. (i)$\Leftrightarrow$(iii) is given by Theorem 4.5 of \cite{SLY16}. $\hfill \Box$

By examining the variation of $ {J}_{\rm soc}^{\rm F}$, we obtain some necessary and sufficient conditions for
the existence of centralized optimal control of {(P1)}.

\begin{theorem}\label{thm1}
	%Suppose that Problem (P1$^{\prime\prime}$) has a minimizer $\check{f}$. Then the following equation systems admits a (unique) solution $(x_i,p_i)$:
	%Suppose  remark or footnote
Assume A1) holds. Then we have the following results:

(i) Problem (P1) has a
%unique
set of optimal control laws %in ${\cal U}_{c}$
if and only if Problem (P1) is convex in $u$ and
	the following equation system admits a set of solutions $(x_i,p_i, \beta_{i}^{j},i,j=1,\cdots,N)$:
	\begin{equation}\label{eq4a}
	\left\{
	\begin{aligned}
	dx_i= &\big(Ax_i+B{u}_i+Gx^{(N)}+f\big)dt+(Cx_i+Du_i+\sigma) dW_i,\\
	dp_i= &-\big(A^Tp_i+C^T\beta_i^i+G^Tp^{(N)}\big)dt-\big(Qx_i-Q_{\Gamma}x^{(N)}-\bar{\eta} \big)dt
	+\sum_{j=1}^N\beta_i^jdW_j,\\
	x_i(0)&={x_{i0}},\ p_i(T)=Hx_i-H_{\Gamma_0}x^{(N)}-\bar{\eta}_0,\ 1\leq i\leq N,
	\end{aligned}\right.
	\end{equation}
	where $p^{(N)}=\frac{1}{N}\sum_{i=1}^Np_i$, and the optimal control $u_i, 1\leq i\leq N$ satisfies the stationary condition
\begin{equation}\label{eq7a}
  R{u}_i+B^Tp_i+D^T\beta_i^i=0.
  \end{equation}

(ii) If Problem (P1) is uniformly convex, then (P1) admits a
%unique
set of optimal control laws.
\end{theorem}
{\it Proof.} See Appendix \ref{app1}. $\hfill{\Box}$%The proof of (i) is similar to Theorem 1 in \cite{W18}; due to page limitation, we omit it.

%Let $\check{x}^{(N)}=\frac{1}{N}\sum_{i=1}^N\check{x}_i$, $\check{u}^{(N)}=\frac{1}{N}\sum_{i=1}^N\check{u}_i$,  and $\check{p}^{(N)}=\frac{1}{N}\sum_{i=1}^N\check{p}_i$.
To ensure the solvability of the problem (P1), we assume

\textbf{A2)} Problem (P1) is uniformly convex in $u$.

We now use the idea inspired by \cite{ZX17}, \cite{ZQ16} to solve the FBSDE (\ref{eq4a}).
Let $p_i=P_Nx_i+K_Nx^{(N)}+s_N$. It follows from (\ref{eq4a}) that
\begin{equation}\label{eq10}
\left\{
\begin{aligned}
d{x}^{(N)}=&\big[(A+G){x}^{(N)}+Bu^{(N)}+f\big]dt+\frac{1}{N}\sum_{i=1}^N(Cx_i+Du_i+\sigma) dW_i,\\
d{p}^{(N)}=&-\Big[(A+G)^T{p}^{(N)}+\frac{1}{N}\sum_{i=1}^NC^T\beta_i^i+(I-\Gamma)^TQ(I-\Gamma){x}^{(N)}
-\bar{\eta}\Big]dt\\
&+\frac{1}{N}\sum_{i=1}^N\sum_{j=1}^N{\beta}_i^jdW_j,\\
{x}^{(N)}(0)=&\frac{1}{N}\sum_{i=1}^Nx_{i0},\ {p}^{(N)}(T)=(H-H_{\Gamma_0})x^{(N)}-\bar{\eta}_0.
\end{aligned}\right.
\end{equation} Then by (\ref{eq4a}), (\ref{eq10}) and It\^{o}'s formula,
\begin{align}\label{s12}
dp_i=&\dot{P}_N{x}_i+P_N\big[\big(A{x}_i+Bu_i+Gx^{(N)}+f\big)dt
+(Cx_i+Du_i+\sigma) dW_i\big]+(\dot{s}_N+\dot{K}_Nx^{(N)})dt\nonumber\\
&+K_N\Big\{\big[(A+G){x}^{(N)}+Bu^{(N)}+f\big]dt
+\frac{1}{N}\sum_{i=1}^N(Cx_j+Du_j+\sigma) dW_j\Big\}\nonumber\\
=&-\big[A^T(P_Nx_i+K_Nx^{(N)}+s_N)
+G^T((P_N+K_N)x^{(N)}+s_N)
+C^T\beta_i^i\nonumber\\&+Qx_i-Q_{\Gamma}x^{(N)}-\bar{\eta} \big]dt+\sum_{j=1}^N\beta_i^jdW_j.
\end{align}
This implies that $\beta_i^i=(P_N+\frac{1}{N}K_N)(Cx_i+Du_i+\sigma)$, and $\beta_i^j=\frac{1}{N}K_N(Cx_j+Du_j+\sigma), \ j\not=i$.
From the stationary condition (\ref{eq7a}),
\begin{equation}\label{eq11a}
Ru_i +B^T(P_Nx_i+K_Nx^{(N)}+s_N)
+D^T(P_N+\frac{1}{N}K_N)(Cx_i+Du_i+\sigma)=0.
\end{equation}
Let $\Upsilon_N\stackrel{\Delta}{=}R+D^T\big(P_N+\frac{K_N}{N}\big)D$. If (\ref{eq11a}) admits a solution, then the optimal control can be given by
\begin{align}\label{eq10e}
 u_i=&-\Upsilon_N^{\dag}\Big[\Big(B^TP_N+D^T\big(P_N+\frac{K_N}{N}\big)C\Big)x_i+B^TK_Nx^{(N)}
 +B^Ts_N+D^T\big(P_N+\frac{K_N}{N}\big)\sigma\Big].%+(I-\Upsilon_N^{\dag}\Upsilon_N)(\Lambda_Nx_i+\Theta_Nx^{(N)}+v_N),
\end{align}
%where $\Lambda_N, \Theta_N\in \mathbb{R}^{n\times n}$ and $v_N\in L^2([0, T], \mathbb{R}^{n})$  are arbitrary.
This together with (\ref{s12}) gives
\begin{align}\label{eq8a}
&\dot{P}_N+A^TP_N+P_NA+C^T(P_N+\frac{K_N}{N})C+Q-\Big(B^TP_N+D^T\big(P_N+\frac{K_N}{N}\big)C\Big)^T\nonumber\\
&\times\Upsilon_N^{\dag}
\Big(B^TP_N+D^T\big(P_N+\frac{K_N}{N}\big)C\Big)=0,\ P_N(T)=H,\\
\label{eq9a}
&\dot{K}_N+(A+G)^TK_N+K_N(A+G)-K_NB\Upsilon_N^{\dag}B^TK_N\nonumber\\
&-\Big(B^TP_N+D^T\big(P_N+\frac{K_N}{N}\big)C\Big)^T\Upsilon_N^{\dag}B^TK_N+G^TP_N+P_NG\nonumber\\
&-K_NB\Upsilon_N^{\dag}\Big(B^TP_N+D^T\big(P_N+\frac{K_N}{N}\big)C\Big)
-Q_{\Gamma}=0,\ K_N(T)=-H_{\Gamma_0},\\
\label{eq10a}
&\dot{s}_N+\Big[A+G-B\Upsilon_N^{\dag}\Big(B^T(P_N+K_N)
+D^T\big(P_N+\frac{K_N}{N}\big)C\Big)\Big]^Ts_N+(P_N+K_N)f-\bar{\eta}\nonumber\\
&+\Big[C-D\Upsilon_N^{\dag}\Big(B(P_N+K_N+D^T\big(P_N+\frac{K_N}{N}\big)C\Big)\Big]^T
\big(P_N+\frac{1}{N}K_N\big)\sigma=0,s_N(T)=-\bar{\eta}_0.
\end{align}

From the above discussion combined with Theorem \ref{thm1}, we have the following result.
\begin{proposition}\label{thm2}
	Assume that A1)-A2) hold. If (\ref{eq8a})-(\ref{eq10a}) admit solutions such that
$$\begin{aligned}
&\mathcal{R}\Big(B^TP_N+D^T\big(P_N+\frac{K_N}{N}\big)C\Big)\cup\mathcal{R}(B^TK_N)\subseteq\mathcal{R}(\Upsilon_N),\\
&B^Ts_N+D^T\big(P_N+\frac{K_N}{N}\big)\sigma\in \mathcal{R}(\Upsilon_N),\ \Upsilon_N\geq0,
\end{aligned}$$
  then Problem (P1) has an optimal control given by (\ref{eq10e}). %$$u_i=-\Upsilon_N^{-1}\left[\Big(B^TP_N+D^T\big(P_N+\frac{K_N}{N}\big)C\Big)x_i+B^TK_Nx^{(N)}+B^Ts_N+D^T\big(P_N+\frac{K_N}{N}\big)\sigma\right].$$
%	where $P_N, K_N $ and $ s _N$ are determined by (\ref{eq8a})-(\ref{eq10a}).
\end{proposition}

%\begin{theorem}\label{thm2}
%	Assume that A1) holds and $Q\geq0$. Then Problem (P1) has an optimal control $$\check{u}_i=-\Upsilon_N^{-1}\left[\Big(B^TP_N+D^T\big(P_N+\frac{K_N}{N}\big)C\Big)x_i+B^TK_Nx^{(N)}+B^Ts_N+D^T\big(P_N+\frac{K_N}{N}\big)\sigma\right],$$
%	where $P_N, K_N $ and $s_N$ are determined by (\ref{eq8a})-(\ref{eq10a}).
%\end{theorem}
%\emph{Proof.} Denote $\Pi=P_N+K_N$. Then from (\ref{eq9a}) and (\ref{eq10a}), $\Pi$ satisfies
%\begin{equation}\label{eq11}
%\begin{aligned}
%\dot{\Pi}&+(A+G)^T\Pi+\Pi (A+G)-\Pi BR^{-1}B^T\Pi+(I-\Gamma)^TQ(I-\Gamma)=0, \quad \Pi(T)=0.
%\end{aligned}
%\end{equation}
%Note that $Q\geq0$ and $R>0$. By \cite{AM90, YZ99},  (\ref{eq8a}) and (\ref{eq11}) admit unique solutions $P\geq 0$
%and $\Pi\geq 0$, respectively, which implies that (\ref{eq9a}) and (\ref{eq10a}) have unique solutions $K_N$ and $s_N$, respectively.
%Then by \cite{MY99, ZX17}, (\ref{eq10}) admits a unique solution. By virtue of this solution, the FBSDE (\ref{eq4a}) is decoupled and the existence %and uniqueness
%of a solution follows.
%By Theorem \ref{thm1}, Problem (P1) has an optimal control
%%the centralized optimal control is
%given by
%$\check{u}_i=-{R^{-1}}B^T(Px_i+Kx^{(N)}+s_N),$
%where $P, K $ and $s_N$ are determined by (\ref{eq8a})-(\ref{eq10a}).  \rightline{$\Box$}

Let $P,K,s$ satisfy
\begin{align}\label{eq8b}
&\dot{P}+A^TP+PA+C^TPC+Q\nonumber
\\
&-\big(B^TP+D^TPC\big)^T
\Upsilon^{\dag}\big(B^TP+D^TPC\big)=0,\quad P(T)=H,\\
\label{eq9b}
&\dot{K}+(A+G)^TK+K(A+G)+G^TP+PG-(B^TP+D^TPC)^T\Upsilon^{\dag}B^TK\nonumber\\
&-KB\Upsilon^{\dag}(B^TP+D^TPC)
-KB\Upsilon^{\dag}B^TK-Q_{\Gamma}=0,\quad K(T)=-H_{\Gamma_0},\\
\label{eq10b}
&\dot{s}+\big[A+G-B\Upsilon^{\dag}\big(B^T(P+K)+D^TPC\big)\big]^Ts+(P+K)f\nonumber\\
&+\big[C-D\Upsilon^{\dag}\big(B(P+K)+D^TPC\big)\big]^TP\sigma
-\bar{\eta}=0,\quad s(T)=-\bar{\eta}_0,
\end{align}
where $\Upsilon\stackrel{\Delta}{=}R+D^TPD$.
For further analysis, we assume

\textbf{{A3)}} (\ref{eq8b})-(\ref{eq10b}) have solutions such that $\Upsilon\geq0$, and
\begin{equation}\label{eq19}
  \mathcal{R}\Big(B^TP+D^TPC\Big)\cup\mathcal{R}(B^TK)\subseteq\mathcal{R}(\Upsilon),\
  B^Ts+D^TP\sigma\in \mathcal{R}(\Upsilon).
\end{equation}
%and $\Upsilon_N\geq0$.
%\footnote{When $\Upsilon_N>0,\Upsilon_N^{\dag}=\Upsilon_N^{-1}. In [], we only consider the case  }
\begin{remark}
{If (\ref{eq8b})-(\ref{eq10b}) have solutions such that $\Upsilon>0$, then $ \Upsilon^{\dag}=\Upsilon^{-1}$ and $\mathcal{R}(\Upsilon)=\mathbb{R}^n$. Thus, assumption A3) holds necessarily. This
corresponds to the case considered in \cite{WZ19}.  }
\end{remark}

As an approximation to ${x}^{(N)}$ in (\ref{eq10}), we obtain
\begin{equation}\label{eq12a}
\frac{d\bar{x}}{dt}=(A+G)\bar{x}-B{\Upsilon^{\dag}}[B^T(P+K)+D^TPC]\bar{x}
-B{\Upsilon^{\dag}}(B^Ts+D^TP\sigma)+f,\ \bar{x}(0)=\bar{x}_0.
\end{equation}
 Then, by Proposition \ref{thm2}, the decentralized control law for agent $i$ can be taken as
\begin{equation}\label{eq15a}
\hat{u}_i(t)=-\Upsilon^{\dag}(t)\big[(B^TP(t)+D^TP(t)C)x_i(t)
+B^TK(t)\bar{x}(t)+B^Ts(t)+D^TP(t)\sigma(t)\big],
%&+
%(I-\Upsilon^{\dag}\Upsilon)(\Lambda x_i+\Theta \bar{x}+v),
\end{equation}
where $P, K$, $s$ and $\bar{x}$ are determined by (\ref{eq8b})-(\ref{eq12a}). %and $\Lambda, \Theta\in \mathbb{R}^{n\times n}$, $v\in L^2([0, T], \mathbb{R}^{n})$ are arbitrary.
\begin{remark}
{ In previous works \cite{HCM12}, \cite{WZ17},
the mean field term $x^{(N)}$ in cost functions (dynamics) is first substituted by a deterministic function $\bar{x}$. By solving an optimal tracking problem subject to consistency requirements, a fixed-point equation of $\bar{x}$ is obtained.
The decentralized control is constructed by handling the fixed-point equation. Here, we first obtain the centralized solution by the variational analysis, and then design decentralized control laws
  by tackling the FBSDEs combined with mean field approximations. Note that in this case $s$ and $\bar{x}$ are fully decoupled and no fixed-point equation is needed.}
\end{remark}

\begin{remark}\label{rem3.2}
  By the local Lipschitz continuous property of the quadratic function, (\ref{eq8b})-(\ref{eq9b}) must admit a unique local solution in a small time duration $[T_0, T]$.
The global existence of the solution for $t\in [0, T]$ can be referred to \cite{AFIJ03}. Particularly, if $Q\geq0$ and $R> 0$, then (\ref{eq8b})-(\ref{eq9b}) admits solutions such that $\Upsilon>0$. Indeed, letting $\Pi=P+K$, $\Pi$ satisfies the following equation
\begin{align}\label{eq11}
&\dot{\Pi}+(A+G)^T\Pi+\Pi (A+G)\
-\big(B^T\Pi+D^TPC\big)^T
\Upsilon^{\dag}\big(B^T\Pi+D^TPC\big)\cr
& +C^TPC+(I-\Gamma)^TQ(I-\Gamma)=0,\quad \Pi(T)=0.
\end{align}
By \cite{YZ99}, if $Q\geq0$ and $R> 0$, then (\ref{eq8b}) and (\ref{eq11}) admit solutions such that $\Upsilon>0$, which implies (\ref{eq8b})-(\ref{eq9b}) admit a solution, respectively.
Besides, from \cite{SLY16}, the solvability of (\ref{eq8b})-(\ref{eq9b}) is equivalent to the uniform convexity of two optimal control problems.
\end{remark}

 After the decentralized control laws (\ref{eq15a}) is applied,
 we have the following closed-loop system
\begin{align}\label{eq20}
d\hat{x}_i=&\big[\bar{A}\hat{x}_i-B\Upsilon^{\dag}(B^T(K\bar{x}+s)+D^TP\sigma)+G\hat{x}^{(N)}+f\big]dt\nonumber\\
&+[\bar{C}\hat{x}_i-D\Upsilon^{\dag}(B^T(K\bar{x}+s)+D^TP\sigma)+\sigma] dW_i,
\end{align}
where $\bar{A}\stackrel{\Delta}=A-B\Upsilon^{\dag}(B^TP+D^TPC)$, and $\bar{C}\stackrel{\Delta}=C-D\Upsilon^{\dag}(B^TP+D^TPC)$.

%After the control law (\ref{eq15a}) is applied,
%the closed-loop state $\hat{x}_i$ satisfies
%\begin{align}\label{eq20i}
%d\hat{x}_i=&\big[\bar{A}\hat{x}_i+G\hat{x}^{(N)}+\bar{f}\big]dt+[\bar{C}\hat{x}_i+\bar{\sigma}] dW_i,
%\end{align}
%where $\bar{A}=A-B\Upsilon^{-1}(B^TP+D^TPC)$, $\bar{C}=C-D\Upsilon^{-1}(B^TP+D^TPC)$, $\bar{f}=f-B\Upsilon^{-1}(B^T(K\bar{x}+s)+D^TP\sigma)$ and $\bar{\sigma}=\sigma-D\Upsilon^{-1}(B^T(K\bar{x}+s)+D^TP\sigma)$.

%Denote
%$${\cal U}_{c} =\Big\{(u_1,\cdots,u_N) \big|\
%u_i(t)\in \sigma\{\bigcup_{i=1}^N{\mathcal F}_t^{i}\}\Big\}.$$

\begin{theorem}\label{thm3}
	Let A1)-A3) hold. %Assume that (\ref{eq8b})-(\ref{eq9b}) have solutions such that $\Upsilon>0$. %in $L^2([0,T], \mathbb{R}^{n\times n})$. %For Problem (P1),
Then for Problem (P1), the set of decentralized control laws
	$\{\hat{u}_1,\cdots,\hat{u}_N\}$ given by (\ref{eq15a}) has asymptotic social optimality, i.e.,
	$$\Big|\frac{1}{N}J^{\rm F}_{\rm soc}(\hat{u})-\frac{1}{N}\inf_{u\in L^2_{{\mathcal F}}(0, T; \mathbb{R}^{nr})}J^{\rm F}_{\rm soc}(u)\Big|=O(\frac{1}{\sqrt{N}}).$$
	
\end{theorem}

\emph{Proof.} See Appendix \ref{app2}. $\hfill{\Box}$

\begin{remark}
The works \cite{HCM12}, \cite{WNZ19} considered the above mean field model with positive (semi-) definite $Q$ and $R$ by the fixed point approach. %The solving process is very complicated.
To achieve asymptotic optimality, an additional condition is needed, like well-posedness of a fixed point equation, which is not easy to verify. Note that in the case $Q\geq 0$ and $R>0$, by {Proposition}
  \ref{prop2} and Remark \ref{rem3.2}, assumptions A1)-A3) hold necessarily. Hence, we get rid of the fixed point condition thoroughly.
% 结果非常复杂， 需要两个额外的条件包括算子小于1 很难验证
\end{remark}

%\begin{theorem}
%  we obtain
%$$\begin{aligned}
%  \inf_{u\in {\cal U}_c}J_{soc}(u)&=\sum_{i=1}^N\mathbb{E}\Big\{[x_{i0}-x^{(N)}(0)]^TP[x_{i0}-x^{(N)}(0)]\cr
%  &+[x^{(N)}(0)]^T\Pi x^{(N)}(0)+2s(0)x^{(N)}(0)\Big\}+q(0),
%  \end{aligned}$$
%where \begin{align*}
%  q(0)=&\int_0^{\infty}e^{-\rho t}\big[(N-1)tr(P\sigma (t)\sigma^T(t))+tr(\Pi \sigma(t) \sigma^T(t))\cr&-N\|B^Ts(t)\|^2_{R^{-1}}+2Ns^T(t)f(t)\big]dt.
%  \end{align*}
%\end{theorem}

\section{Mean Field LQ Social Control over an Infinite Horizon}

Based on the similar discussion and analysis in Section \ref{sec3}, we may design the following decentralized control laws
for Problem (P0):
%\footnote{To avoid to verify the condition (\ref{eq19}) in the infinite horizon, we only consider the case $\Upsilon>0$.  }:
\begin{equation}\label{eq14}
\begin{aligned}
\hat{u}_i(t)= &-{\Upsilon^{\dag}}\left[(B^TP+D^TPC)x_i(t)+B^T(\Pi-P)\bar{x}(t)\right.\cr
&\left.+B^Ts(t)+D^TP\sigma(t)\right],\ \ i=1,\cdots, N,
\end{aligned}
\end{equation}
where $\Upsilon=R+D^TPD$, $P$ and $\Pi$ are determined by
\begin{equation}\label{eq15}
  A^TP+PA+C^TPC-\big(B^TP+D^TPC\big)^T
  \Upsilon^{\dag}\big(B^TP+D^TPC\big)+Q=0,
\end{equation}
\begin{equation}\label{eq16}
  (A+G)^T\Pi+\Pi (A+G)-\big(B^T\Pi+D^TPC\big)^T
  \Upsilon^{\dag}\big(B^T\Pi+D^TPC\big)+C^TPC+{Q}-Q_{\Gamma}=0,
\end{equation}
and
$s, \bar{x}\in L_2([0,\infty),\mathbb{R}^n)$ are determined by
\begin{align}\label{eq17}
\ &\frac{ds}{dt}+\big[A+G-B\Upsilon^{\dag}\big(B^T\Pi+D^TPC\big)\big]^Ts+\Pi f
 +\big[C-D\Upsilon^{\dag}\big(B\Pi+D^TPC\big)\big]^TP\sigma-\bar{\eta}=0,\\
 \label{eq18}
&\frac{d\bar{x}}{dt}=\big[A+G-B{\Upsilon^{\dag}}(B^T\Pi+D^TPC)\big]\bar{x}
-B{\Upsilon^{\dag}}(B^Ts+D^TP\sigma)+f,\ \bar{x}(0)=\bar{x}_0.
\end{align}
Here the existence conditions of $P, \Pi,s$ and $\bar{x}$ %need to be investigated further.
 are to be ensured later.

For further analysis, we first introduce some definitions. Consider the following system
\begin{align}\label{eq37a}
  dy(t)&=(Ay(t)+Bu(t))dt+(Cy(t)+Du(t))dW(t),\\
  \label{eq37b}
  z(t)&=Fy(t),
\end{align}
where $y(t)\in \mathbb{R}^n$, and $W(t)$ is a 1-dimensional Brownian motion.
\begin{definition}
The system (\ref{eq37a}) with $u=0$ (or simply $[A,C]$) is said to be mean-square stable, %in the sense,
if
for any initial value $y(0)$, $\lim_{t\to\infty}\mathbb{E}[y^T(t)y(t)]=0$.
\end{definition}

\begin{definition}
The system (\ref{eq37a}) (or simply $[A,B;C,D]$) is said to be stabilizable (in the mean-square sense), if
  there exists a control law $u(t)=Ky(t)$ such that
for any initial $y(0)\in \mathbb{R}^n$, the closed-loop system
$dy(t)=(A+BK)y(t)dt+(C+DK)y(t)dW(t)$
is mean-square stable. In this case $u(t)$ is called a stabilizer.
If $C=D=0$, then the system, abbreviated as $(A,B)$, is stabilizable.
\end{definition}
\begin{definition}\cite{ZZC08} The system (\ref{eq37a})-(\ref{eq37b}) (or simply $[A,C;F]$) is said to be exactly observable, if there exists a $T_0\geq0$ such that
for any $T>T_0$, $z(t)=0, u(t)=0,a.s.,\ 0\leq t\leq T$ implies $y(0)=0$. If $C=0$, then the system,
abbreviated as $(A,F)$, is observable.
\end{definition}

\begin{definition}\cite{ZZC08} The system (\ref{eq37a})-(\ref{eq37b}) (or simply $[A,C;F]$) is said to be exactly detectable, if
there exists a $T_0\geq0$ such that
for any $T>T_0$, $z(t)=0, u(t)=0,a.s.,\ 0\leq t\leq T$ implies $\lim_{t\to\infty}\mathbb{E}[y^T(t)y(t)]=0$.
% If $C=0$, then the system,
%abbreviated as $(A,C)$, is detectable.
\end{definition}

%For the system (\ref{eq37a}), define the cost function
%\begin{equation}\label{eq38}
%  J(u)=\mathbb{E}\int_0^{\infty}(\|y(t)\|_Q^2+\|u(t)\|_R^2)dt.
%  \end{equation}

Some basic assumptions are listed for reference:

\textbf{A4)} The system $[A, B; C,D]$ is stabilizable, and the system $(A+G, B)$ is stabilizable.

\textbf{A5}) %The problem (\ref{eq37a}), (\ref{eq38}) is uniformly convex, % or $T(P)>0$ 用线性矩阵不等式来刻画。
%$\{\bar{P}: \mathcal{H}(\bar{P})\geq0, R+D^T\bar{P}D\geq0\} has nonempty interior
$\mathcal{S}_1=\big\{\bar{P}=\bar{P}^T: \mathcal{H}(\bar{P})\geq0,   \hbox{ker}(R_{\bar{P}})\subseteq  \hbox{ker}(B)\cap  \hbox{ker}(D),[A,C, Q_{\bar{P}}^{1/2}]\hbox{ is exactly detectable}\big\}\\\not=\emptyset$,
$\mathcal{S}_2=\big\{\bar{\Pi}=\bar{\Pi}^T: \mathcal{M}(\bar{\Pi})\geq0,[A+G, Q_{\bar{\Pi}}^{1/2}] \hbox{ is  detectable}\big\}\not=\emptyset$, where
$$\mathcal{H}(\bar{P})=\left[\begin{array}{cc}
Q_{\bar{P}}&\bar{P}B+C^T\bar{P}D \\
B^T\bar{P}+D^T\bar{P}C&R_{\bar{P}}
\end{array}\right],$$
$$\mathcal{M}(\bar{\Pi})=\left[\begin{array}{cc}
Q_{\bar{\Pi}}& \bar{\Pi}B+C ^T PD\\
B^T\bar{\Pi}+ D^TPC& R_{\bar{P}}
\end{array}\right],$$
%$$M=\left[\begin{array}{cc}
%A+G -B\Upsilon^{\dag}D^TPC& B\Upsilon^{\dag}B^T \\
%Q-Q_{\Gamma}+C^TPC-(D^TPC)^T\Upsilon^{-1}D^TPC& -(A+G-B\Upsilon^{-1}D^TPC)^{T}
%\end{array}\right].$$
with
\begin{align*}
  Q_{\bar{P}}&=A^T\bar{P}+\bar{P}A+C^T\bar{P}C+Q,\\
R_{\bar{P}}&=R+D^T\bar{P}D,\\
Q_{\bar{\Pi}}&=(A+G)^T\bar{\Pi}+\bar{\Pi}(A+G)+C^TPC+Q-Q_{\Gamma}.
\end{align*}

%\begin{remark}
%   $M$ is a Hamiltonian matrix. The Hamiltonian matrix plays a significant role in studying %general
%   algebraic Riccati equations. See more details of Hamiltonian matrices in \cite{AFIJ03, M77}.
%\end{remark}

\begin{lemma}\label{lem2a}
	Under A4)-A5), the following holds:

(i) (\ref{eq15}) admits a unique solution $P$ such that $\Upsilon\geq 0$ and $[\bar{A},\bar{C}]$ is mean-square stable, where $\bar{A}=A-B\Upsilon^{\dag}(B^TP+D^TPC)$, and $\bar{C}=C-D\Upsilon^{\dag}(B^TP+D^TPC)$; %R+D^TPD

(ii) (\ref{eq16}) admits a unique solution $ \Pi$ such that $A+G-B{\Upsilon^{\dag}}(B^T\Pi+D^TPC)$ is Hurwitz;

(iii) (\ref{eq17})-(\ref{eq18}) admits a set of unique solutions $s, \bar{x}\in L_2([0,\infty),\mathbb{R}^n)$;

(iv) %The following holds:
$
\mathcal{R}\Big(B^TP+D^TPC\Big)\cup\mathcal{R}(B^T(\Pi-P))\subseteq\mathcal{R}(\Upsilon), B^Ts+D^TP\sigma\in \mathcal{R}(\Upsilon).$
\end{lemma}

\emph{Proof.} Applying Theorem 2 in \cite{LQZ19}, we obtain that under A4)-A5), (\ref{eq15}) admits a unique solution $P$ such that the system
$[\bar{A}, \bar{C}]$ is mean-square stable. Note that in  (\ref{eq16}), $P$ is known. Since $(A+G,B)$ is stabilizable, then from \cite[Theorem 2]{LQZ19}, (\ref{eq16}) admits a unique solution $ \Pi$ such that $A+G-B{\Upsilon^{\dag}}(B^T\Pi+D^TPC)$ is Hurwitz. From an argument in \cite[Appendix A]{WZ12}, we obtain
$s\in L_2([0,\infty),\mathbb{R}^n)$ if and only if
\begin{equation*}
  s(0)=\int_0^{\infty}e^{[A+G-B{\Upsilon^{\dag}}(B^T\Pi+D^TPC)]\tau}
  (\Pi f+\bar{C}P\sigma-\bar{\eta})d\tau.
\end{equation*}
Under this initial condition, we have
\begin{equation*}
  s(t)=\int_t^{\infty}e^{-[A+G-B{\Upsilon^{\dag}}(B^T\Pi+D^TPC)](t-\tau)}
 (\Pi f+\bar{C}P\sigma-\bar{\eta})d\tau.
\end{equation*}
%Since $A+G-B{\Upsilon^{\dag}}(B^T\Pi+D^TPC)$ is Hurwitz, then $\bar{x}\in L_2([0,\infty),\mathbb{R}^n)$.
%Since $\Upsilon\geq 0$, there exists an orthogonal matrix $V$ such that
%$$\Upsilon=V
%\left[\begin{array}{cc}
%\Lambda&0\\
%0&0
%\end{array}\right]V^T,$$
%which implies $$\Upsilon^{\dag}=V
%\left[\begin{array}{cc}
%\Lambda^{-1}&0\\
%0&0
%\end{array}\right]V^T.$$
From the argument in \cite[Theorem 1]{LQZ19}, one can show that $(B^Ts+D^TP\sigma)^T(I-\Upsilon\Upsilon^{\dag})=0$,
which implies $B^Ts+D^TP\sigma\in \mathcal{R} (\Upsilon)$. Similarly, we have $
\mathcal{R}\Big(B^TP+D^TPC\Big)\cup\mathcal{R}(B^TK)\subseteq\mathcal{R}(\Upsilon)$.
$\hfill \Box$

%After the control law (\ref{eq14}) is applied,
%the closed-loop state $\hat{x}_i$ satisfies
%\begin{align}\label{eq20i}
%d\hat{x}_i=&\big[\bar{A}\hat{x}_i-B\Upsilon^{-1}(B^T(K\bar{x}+s)+D^TP\sigma)+G\hat{x}^{(N)}+\bar{f}\big]dt+[\bar{C}\hat{x}_i+\bar{\sigma}] dW_i,
%\end{align}
%where $\bar{f}=f-B\Upsilon^{-1}(B^T(K\bar{x}+s)+D^TP\sigma)$ and $\bar{\sigma}=\sigma-D\Upsilon^{-1}(B^T(K\bar{x}+s)+D^TP\sigma)$.

We now introduce an additional assumption. Later, the assumption is shown to be necessary and sufficient for the uniform stabilization of all the subsystems.

\textbf{A6)} $\bar{A}+G$ is Hurwitz, where $\bar{A}{=}A-B\Upsilon^{\dag}(B^TP+D^TPC)$.

It is shown that the decentralized control laws (\ref{eq15a}) uniformly stabilize the systems (\ref{eq1}) .

\begin{theorem}\label{thm4}
	Let A1), A4)-A6) hold.  Then there exists an $N_0$ such that for $N\geq N_0$, the following hold:
	\begin{align}\label{eq13b}
	%\sup_{N\geq 1}
	&\max_{1\leq i\leq N}\mathbb{E}\int_0^{\infty} \left(\|\hat{x}_i(t)\|^2+\|\hat{u}_i(t)\|^2\right)dt<\infty.\\
&\mathbb{E}\int_0^{\infty} \|\hat{x}^{(N)}(t)-\bar{x}(t)\|^2dt=O(\frac{1}{N}).	\label{eq13c}
\end{align}

\end{theorem}

\emph{Proof.}
After the control (\ref{eq14}) is applied, we have
\begin{align}\label{eq20i}
d\hat{x}_i=&\big[\bar{A}\hat{x}_i+G\hat{x}^{(N)}+\bar{f}\big]dt+[\bar{C}\hat{x}_i+\bar{\sigma}] dW_i,
\end{align}
 where $\bar{f}\stackrel{\Delta}{=}f-B\Upsilon^{\dag}(B^T(K\bar{x}+s)+D^TP\sigma)$, and $\bar{\sigma}\stackrel{\Delta}{=}\sigma-D\Upsilon^{\dag}(B^T(K\bar{x}+s)+D^TP\sigma)$.
 Let $\xi(t)=\hat{x}^{(N)}(t)-\bar{x}(t)$.
From (\ref{eq20i}) and (\ref{eq18}),
\begin{equation}\label{eq21d}
  \xi(t)=e^{(\bar{A}+G)t}\xi(0)+\frac{1}{N}\sum_{i=1}^N\int_0^te^{(\bar{A}+G)(t-\tau)}(\bar{C}\hat{x}_i+\bar{\sigma}) dW_i.
\end{equation}
Thus, we have
\begin{align}\label{eq43e}
&\mathbb{E}\int_0^{T}  \left(\|\hat{x}^{(N)}(t)-\bar{x}(t)\|^2\right)dt\cr
\leq\ & 2\mathbb{E}\int_0^{T}\left\| e^{(\bar{A}+G)t}\right\|^2\big\|\hat{x}^{(N)}(0)-\bar{x}(0)\big\|^2dt
+\!2\mathbb{E}\!\int_0^{T}\!\frac{1}{N}\left\| \int_0^te^{(\bar{A}+G)(t-\tau)}(\bar{C}\hat{x}_i+\bar{\sigma}) dW_i(\tau)\right\|^2dt\cr
\leq\ & 2\int_0^{T} \left\| e^{(\bar{A}+G)t}\right\|^2\mathbb{E}\big\|\hat{x}^{(N)}(0)-\bar{x}(0)\big\|^2dt\
+\frac{2}{N} \mathbb{E}\int_0^{T}\int_0^t \left\| e^{(\bar{A}+G)(t-\tau)}\right\|^2\|\bar{C}\hat{x}_i+\bar{\sigma}\|^2 d\tau dt\cr
\leq\ & \frac{2}{N}\int_0^{T} \left\| e^{(\bar{A}+G)t}\right\|^2\max_{1\leq i\leq N}\mathbb{E}\big\|\hat{x}_i(0)\big\|^2dt
+ \frac{c}{N} \mathbb{E}\int_0^{T}(c_1\|\hat{x}_i\|^2+c_2)\int_\tau^{T}\big\|e^{(\bar{A}+{G})(t-\tau)}\big\|^2 dtd\tau\cr
\leq& \frac{c_1}{N}\max_{1\leq i\leq N}\mathbb{E}\int_0^{T}\|\hat{x}_i\|^2dt+\frac{c_1}{N}.
\end{align}
%By Lemma \ref{lem2a} (iii) and Lemma \ref{lem2},
%\begin{equation}
%  \mathbb{E}\int_0^{\infty} \|\hat{x}^{(N)}\|^2dt\leq  \mathbb{E}\int_0^{\infty} (\|\hat{x}^{(N)}-\bar{x}\|^2+\bar{x})dt<\infty.
%  \end{equation}
Let $P$ satisfy
$$P\bar{A}+\bar{A}^TP+\bar{C}^TP\bar{C}=-2I.$$
%where $Q<0$.
From Lemma \ref{lem2a}(i) and \cite{RCMZ01}, we have $P>0$.
By It\^{o}'s formula and (\ref{eq20i}),
\begin{align}\label{eq40a}
  &\mathbb{E}[\hat{x}_i^T(T)P\hat{x}_i(T)-\hat{x}_i^T(0)P\hat{x}_i(0)]\cr
  =&\mathbb{E}\int_0^T\big[\hat{x}_i^TP(\bar{A}\hat{x}_i+G\hat{x}^{(N)}+\bar{f})
  +(\bar{A}\hat{x}_i+G\hat{x}^{(N)}+\bar{f})^TP\hat{x}_i\big]dt\cr
  &+\mathbb{E}\int_0^T(\bar{C}\hat{x}_i+\bar{\sigma})^TP(\bar{C}\hat{x}_i+\bar{\sigma})dt.
\end{align}
From (\ref{eq40a}), we have
\begin{align}\label{eq43b}
  &\mathbb{E}[\hat{x}_i^T(T)P\hat{x}_i(T)-\hat{x}_i^T(0)P\hat{x}_i(0)]\cr
  =&\mathbb{E}\int_0^T\big[\hat{x}_i^T(P\bar{A}+\bar{A}^TP+\bar{C}^TP\bar{C})\hat{x}_i
  +(\hat{x}^{(N)})^T(PG+G^TP)\hat{x}^{(N)}\cr
  &+2(P\bar{f}+\bar{C}^TP\bar{\sigma})^T\hat{x}_i+\bar{\sigma}^TP\bar{\sigma}\big]dt\cr
  \leq &\mathbb{E}\int_0^T\big[\hat{x}_i^T(P\bar{A}+\bar{A}^TP+\bar{C}^TP\bar{C})\hat{x}_i+\|\hat{x}_i\|^2
  +(\hat{x}^{(N)})^T(PG+G^TP)\hat{x}^{(N)}\cr
  &+\|P\bar{f}+\bar{C}^TP\bar{\sigma}\|^2+ \bar{\sigma}^TP\bar{\sigma}\big]dt\cr
  \leq& -\mathbb{E}\int_0^T(\hat{x}_i^T\hat{x}_i)dt+\alpha_T,
%\leq-\frac{1}{\lambda_{\max}(P)}\sum_{i=1}^N
%\mathbb{E}\int_0^T(\hat{x}_i^TP\hat{x}_i)dt+\alpha_T,
\end{align}
where %$\lambda_{\max}(P)$ is the maximum eigenvalue of $P$ and
\begin{equation*}
  \alpha_T=\mathbb{E}\int_0^T\big[(\hat{x}^{(N)})^T(PG+G^TP)\hat{x}^{(N)}
  +\|P\bar{f}+\bar{C}^TP\bar{\sigma}\|^2+ \bar{\sigma}^TP\bar{\sigma}\big]dt.
\end{equation*}
This with (\ref{eq43e}) gives
\begin{align}
 \mathbb{E}\int_0^T\|\hat{x}_i\|^2dt\leq& \mathbb{E}[{x}_{i0}^TP{x}_{i0}]+\alpha_T\leq c_2\mathbb{E}\int_0^T\|x^{(N)}\|^2dt+c_2\cr
  \leq& 2c_2\mathbb{E}\int_0^T\big(\|\bar{x}(t)\|^2+\|\xi(t)\|^2\big)dt+c_2\cr
  \leq& 2c_2\Big[\mathbb{E}\int_0^T\|\bar{x}(t)\|^2dt+\frac{c_1}{N}\max_{1\leq i\leq N}\mathbb{E}\int_0^{T}\|\hat{x}_i\|^2dt\Big]
  +\frac{2c_1c_2}{N}+c_2.
  \end{align}
Thus, there exists $N_0$ such that for any $N>N_0$,
$$\max_{1\leq i\leq N}\mathbb{E}\int_0^T\|\hat{x}_i\|^2dt\leq 2c_2\mathbb{E}\int_0^T\|\bar{x}(t)\|^2dt+{2c_1c_2}+c_2.$$
  Note $\bar{x}\in  L_2([0,\infty),\mathbb{R}^n)$.
We have $$\max_{1\leq i\leq N}\mathbb{E}\int_0^{\infty}\|\hat{x}_i\|^2dt\leq c. $$
This together with (\ref{eq43e}) gives (\ref{eq13c}).
  \hfill $\Box$

We now give two equivalent conditions for uniform stabilization of all the subsystems.

\begin{theorem}\label{thm5}
	For (P0), let A5) hold. Assume that (\ref{eq15})-(\ref{eq16}) have solutions. Then for (P0) the following statements are equivalent:

 (i) there exists an $N_0$ such that for $N\geq N_0$ and any initial condition $(\hat{x}_1(0),\cdots, \hat{x}_N(0))$ satisfying A1),
	\begin{equation}\label{eq23}
	%\sup_{N\geq 1}
	\sum_{i=1}^N\mathbb{E}\int_0^{\infty} \left(\|\hat{x}_i(t)\|^2+\|\hat{u}_i(t)\|^2\right)dt<\infty;
	\end{equation}	
	(ii) (\ref{eq15})-(\ref{eq17}) admit solutions such that $R+D^TPD\geq0$, $
\mathcal{R}\Big(B^TP+D^TPC\Big)\cup\mathcal{R}(B^T(\Pi-P))\subseteq\mathcal{R}(\Upsilon), B^Ts+D^TP\sigma\in \mathcal{R}(\Upsilon),$
 and $\bar{A}+G$ is Hurwitz;
	
(iii) A4) and A6) hold.
	
\end{theorem}

\emph{Proof.} See Appendix \ref{app3}. $\hfill{\Box}$

For the case $Q\geq0$, $R>0$, when the assumption A5) is strengthened to A5)$^{\prime}$, we can give the following equivalent conditions for uniform stabilization of the systems.

${\bf{A5^{\prime})}}$
$Q\geq0$, $R>0$, $[A,C, \sqrt{Q}]$ is exactly observable, and
$(A+G, \sqrt{Q}(I-\Gamma))$ is observable.

%Note that assumptions A2)-A3) imply $(A+G,B)$ is stabilizable and $(A+G,C^TPC+Q-Q_{\Gamma})$ is observable.

\begin{theorem}\label{thm5b}
	Let {A5$^{\prime}$)} hold. Assume that (\ref{eq15})-(\ref{eq16}) have solutions. Then for (P0) the following statements are equivalent:
	
	(i) For any initial condition $(\hat{x}_1(0),\cdots, \hat{x}_N(0))$ satisfying A1), the following holds,
	\begin{equation*}%\label{eq23}
	%\sup_{N\geq 1}
	\sum_{i=1}^N\mathbb{E}\int_0^{\infty} \left(\|\hat{x}_i(t)\|^2+\|\hat{u}_i(t)\|^2\right)dt<\infty;
	\end{equation*}
	(ii) (\ref{eq15}) and (\ref{eq16}) admit unique solutions such that $P>0, \Pi>0$, and $\bar{A}+G$ is Hurwitz;
	
	(iii) A4) and A6) hold.
	
\end{theorem}

%\section{Proofs of Theorems \ref{thm5} and \ref{thm6}}\label{app c}
%\def\theequation{C.\arabic{equation}}
%\setcounter{equation}{0}

\emph{Proof.} See Appendix \ref{app3}. $\hfill{\Box}$

\begin{remark}
{ In \cite{QZW19}, some similar results were given for the stabilization of mean field systems. However, only the limiting problem
  was considered in their work and the mean field term in dynamics and costs is $\mathbb{E}x(t)$ instead of $x^{(N)}(t)$. Here we study large-population multiagent systems
  and the number of agents is large but not infinite. %Approximation errors for large-population systems are further analyzed.
  The errors of mean field approximations need to be further analyzed.
  In this case, an additional assumption
  A6) is needed to obtain uniform stabilization.}
\end{remark}

To compare the optimal social costs under decentralized and centralized strategies, we need the presumption that Problem (P0)
admits a centralized solution. Thus, we set an assumption on the following generalized Riccati equation:

\textbf{A7)}
The equation
\begin{equation*}
\check{\textbf{A}}^T\textbf{P}+\textbf{P}\check{\textbf{A}}+\sum_{i=1}^N\textbf{C}_i^T\textbf{P}\textbf{C}_i+\bar{\textbf{Q}}-
\Big(\textbf{B}^T\textbf{P}+\sum_{i=1}^N\textbf{D}_i^T\textbf{P}\textbf{C}_i\Big)^T
{\mathbf{\Upsilon}}^{\dag}
\Big(\textbf{B}^T\textbf{P}+\sum_{i=1}^N\textbf{D}_i^T\textbf{P}\textbf{C}_i\Big)=0
\end{equation*}
 admits a solution such that $\mathbf{\Upsilon}=\textbf{R}
+\sum_{i=1}^N\textbf{D}_i^T\textbf{P}\textbf{D}_i\geq0$, $\mathcal{R}\Big(\textbf{B}^T\textbf{P}+\sum_{i=1}^N\textbf{D}_i^T\textbf{P}\textbf{C}_i\Big)\subseteq \mathcal{R}(\mathbf{\Upsilon})$ and the following system is mean-square stable:
\begin{equation*}
  d\textbf{x}=\Big[\check{\textbf{A}}-\textbf{B}\mathbf{\Upsilon}^{\dag}\Big(\textbf{B}^T\textbf{P}+\sum_{i=1}^N\textbf{D}_i^T\textbf{P}\textbf{C}_i\Big)\Big]\textbf{x}dt +\sum_{i=1}^N\Big[\textbf{C}_i-\textbf{D}_i\mathbf{\Upsilon}^{\dag}\Big(\textbf{B}^T\textbf{P}+\sum_{i=1}^N\textbf{D}_i^T\textbf{P}\textbf{C}_i\Big)\Big]dW_i.
\end{equation*}
  %where $\mathbf{\Upsilon}=\textbf{R}+\sum_{i=1}^N\textbf{D}_i^T\textbf{P}\textbf{D}_i$.

%\textbf{(A6)} Problem (PS) is uniformly convex.

We now are in a position to state the asymptotic social optimality of
the decentralized control.
\begin{theorem}\label{thm8}
	Let A1), A4)-A7) hold. For Problem (P0), the set of decentralized control laws
	$\{\hat{u}_1,\cdots,\hat{u}_N\}$ given by (\ref{eq14}) has asymptotic social optimality, i.e.,
	$$\Big|\frac{1}{N}J_{\rm soc}(\hat{u})-\frac{1}{N}\inf_{u\in \mathcal{U}_c}J_{\rm soc}(u)\Big|=O(\frac{1}{\sqrt{N}}).$$
	
\end{theorem}

We first provide a preliminary lemma, which plays an important role in showing asymptotic optimality of decentralized control.
\begin{lemma}\label{lem1.3}
  For the system (\ref{eq37a}), assume $[A,B;C,D]$ is stabilizable. Then for %any $K\in \mathbb{R}^{n\times n}$ and
  any $u\in L_2([0,\infty),\mathbb{R}^n)$ and a stabilizer $Ky$, %, %making sure $y\in L_2([0,\infty),\mathbb{R}^n)$,
  there exist constants $\alpha_i,c_i>0, i=1,2$ such that
  $$\mathbb{E}\int_0^{\infty}\|y(t)\|^2dt\leq \alpha_1 \mathbb{E} \int_0^{\infty}\|u(t)-Ky(t)\|^2dt+c_1,$$
  $$\mathbb{E}\int_0^{\infty}\|u(t)\|^2dt\leq \alpha_2 \mathbb{E} \int_0^{\infty}\|u(t)-Ky(t)\|^2dt+c_2.$$
\end{lemma}

\emph{Proof.} %Suppose $[A,C]$ is mean-square stable.
%Let $Ky\in L_2([0,\infty),\mathbb{R}^n)$ be a stabilizer of the system (\ref{eq37a}).
Define
%$\mathcal{T}: L_2([0,\infty),\mathbb{R}^n)\rightarrow L_2([0,\infty),\mathbb{R}^n)$ by
$u^*=u-Ky$, where $y$ satisfies (\ref{eq37a}). Then $u^*\in L_2([0,\infty),\mathbb{R}^n)$ and
 $y$ satisfies
%By \cite{SY17}, there exists $c$ such that $\mathbb{E}\int_0^{\infty}\|y(t)\|^2dt\leq c\mathbb{E}\int_0^{\infty}\|u(t)\|^2dt.$
%Then the operator $\mathcal{T}$ is bounded. The inverse $\mathcal{T}^{\dag}$ is given by $\mathcal{T}^{\dag}u^*=u^*+Ky^*$,
%where $y^*$ is the solution to the equation
\begin{equation*}
  dy(t)=[(A+BK)y(t)+Bu^*(t)]dt
+[(C+DK)y(t)+u^*(t)]dW(t),\ y(0)=y_0.
\end{equation*}
Since $Ky$ is a stabilizer, then by \cite{SY17}, there exists a constant $\alpha_1$ such that $\mathbb{E}\int_0^{\infty}\|y(t)\|^2dt\leq \alpha_1\mathbb{E}\int_0^{\infty}\|u^*(t)\|^2dt+c_1.$
Hence,
\begin{align*}
  \mathbb{E}\int_0^{\infty}\|u(t)\|^2dt=&  \mathbb{E}\int_0^{\infty}\|u^*(t)+Ky(t)\|^2dt\\
\leq& \alpha_2 \mathbb{E} \int_0^{\infty}\|u^*(t)\|^2dt +c_2\\
=&\alpha_2 \mathbb{E} \int_0^{\infty}\|u(t)-Ky(t)\|^2dt+c_2,
\end{align*}
where $\alpha_2=2\alpha_1\|K\|^2+2,$ and $c_2=2c_1$.
\hfill $\Box$

\emph{Proof of Theorem \ref{thm8}.} We first prove that for $u\in \mathcal{U}_c$, $J_{\rm soc}(u)< c_1$ implies that there exists a constant $c_2$ such that
\begin{equation}\label{eq36a}
\mathbb{E}\int_0^{\infty}(\|x_i\|^2+\|u_i\|^2)dt<c_2,
\end{equation} for all $i=1,\cdots,N$.
From A7), the following equation admits a unique solution $\textbf{s}\in L_2([0,\infty),\mathbb{R}^{Nn})$,
$$\begin{aligned}
\dot{\textbf{s}}&+\big[\check{\textbf{A}}-\sum_{i=1}^N\textbf{B}\mathbf{\Upsilon}^{\dag}\big(\textbf{B}^T\textbf{P}+\textbf{D}_i^T\textbf{P}\textbf{C}_i\big)\big]^T\textbf{s}
+\textbf{P} (f\otimes \textbf{1})\\
&\ +\sum_{i=1}^N\big[\textbf{C}_i-\textbf{D}_i\mathbf{\Upsilon}^{\dag}\big(\textbf{B}^T\textbf{P}+\textbf{D}^T_i\textbf{P}\textbf{C}_i\big)\big]^T\textbf{P}\sigma_i-\bar{\eta}\otimes \textbf{1}=0.
  \end{aligned}$$
By It\^{o}'s formula, we have
\begin{align*}
 J_{\rm soc}(u)&=\limsup_{T\to\infty}\mathbb{E}\big[{\textbf{x}}^T(0)\textbf{P}{\textbf{x}}(0)
 -{\textbf{x}}^T(T)\textbf{P}{\textbf{x}}(T)\big]\cr
 &\quad+\mathbb{E}\int_0^{\infty}\Big\|{\textbf{u}}+
 {\bf{\Upsilon}}^{\dag}\Big[\big(\textbf{B}^T\textbf{P}+\sum_{i=1}^N\textbf{D}_i^T\textbf{P}\textbf{C}_i\big)
 {\textbf{x}}
 \quad+\textbf{B}^T\textbf{s}+\sum_{i=1}^N\textbf{D}^T_i
 \textbf{P}\sigma_i\Big]\Big\|^2dt\cr
 &\geq\mathbb{E}\int_0^{\infty}\Big\|{\textbf{u}}+
 {\bf{\Upsilon}}^{\dag}\big(\textbf{B}^T\textbf{P}+\sum_{i=1}^N\textbf{D}_i^T\textbf{P}\textbf{C}_i\big)
 {\textbf{x}}\Big\|^2dt-c.
  \end{align*}
  By Lemma \ref{lem1.3}, there exist constants $\alpha,c>0$ such that
\begin{align}\label{eq36}
&\sum_{i=1}^N\mathbb{E}\int_0^{\infty}(\|x_i\|^2+\|u_i\|^2)dt\cr
\leq &\alpha\mathbb{E}\int_0^{\infty}\Big\|{\textbf{u}}+
 {\bf{\Upsilon}}^{\dag}(\textbf{B}^T\textbf{P}+\sum_{i=1}^N\textbf{D}_i^T\textbf{P}\textbf{C}_i){\textbf{x}}\Big\|^2dt+c
\leq \alpha J_{\rm soc}(u)+c\leq c_2.
\end{align}
By a similar argument to the proof of Theorem \ref{thm3} combined with Theorem \ref{thm4}, the conclusion follows.
$\hfill \Box$

\begin{remark}

If A5) is replaced by A5$^{\prime}$), then it can be shown that the decentralized control (\ref{eq14}) still has asymptotic social optimality.
% If G=0, C=D=0 我们的模型变成Huang TAC, 我们去掉了不动点假设这个条件。。。
\end{remark}

\section{Asymptotically Social Optimal Cost}
We now give an explicit expression of the asymptotic average social optimum in terms of the solutions of two Riccati equations.

\begin{theorem}\label{thm5.1}
Assume i) A1), A4-A7) hold; ii) $\{x_{i0}\}$ have the same variance. Then the asymptotic average social optimum is given by
\begin{equation*}
  \lim_{N\to\infty}\frac{1}{N}J_{\rm soc}(\hat{u})=\mathbb{E}\big[(x_{i0}-\bar{x}_0)^TP(x_{i0}-\bar{x}_0)
  +\bar{x}_0^T\Pi \bar{x}_0+2s^T(0)\bar{x}_0\big]+m,
\end{equation*}
\end{theorem}
where $P$ and $ \Pi$ are given by (\ref{eq15})-(\ref{eq16}), respectively, and
\begin{equation*}
  m=\int_0^{\infty}\big[\|\sigma(t)\|^2_P-\|B^Ts(t)+D^TP\sigma(t)\|^2_{\Upsilon^{\dag}}
  %-s^TB\Upsilon^{\dag}D^TP\sigma
  +2s^T(t)f(t)+\|\eta(t)\|^2_Q\big]dt.
\end{equation*}
To prove Theorem \ref{thm5.1}, we need two lemmas.

Consider the mean-field type system
\begin{equation}\label{eq63}
  dz_i=(Az_i+Bu_i+G\mathbb{E}[{z}_i]+f)dt
  +(Cz_i+Du_i+\sigma)dW_i,\ z_i(0)=x_{i0},
\end{equation}
with the cost function
\begin{equation}\label{eq63b}
\mathcal{J}_i(u_i)=\mathbb{E}\int_0^{\infty}(\|z_i-\Gamma \mathbb{E}[{z}_i]-\eta\|^2_Q+\|u_i\|^2_R)dt.
\end{equation}
The admissible control set is given by
\begin{equation*}
{\mathcal U}_{i} =\Big\{u_i \big|\
u_i(t)\ \hbox{is adapted to } \sigma(z_i(s),0\leq s\leq t), \mathbb{E}\int_0^{\infty}\|z_i(t)\|^2dt<\infty,\forall i \Big\}.
%u_i(t)\ \hbox{is adapted to}\  \sigma(x_i(s),0\leq s\leq t),  \mathbb{E}\int_0^{\infty}\|x_i(t)\|^2dt<\infty
\end{equation*}
\begin{lemma}\label{lem5.1}
   For the system (\ref{eq63})-(\ref{eq63b}),
the optimal control is given by
\begin{equation}\label{eq64}
 \hat{u}_i=-\Upsilon^{\dag}[(B^TP+D^TPC)z_i
  +B^T(\Pi-P)\mathbb{E}[{z}_i]+B^Ts+D^TP\sigma],
\end{equation}
and the optimal cost is
\begin{equation*}
  \inf_{u_i\in {\mathcal U}_i}\mathcal{J}_{i}(u_i)=\mathbb{E}\big[(x_{i0}-\bar{x}_0)^TP(x_{i0}-\bar{x}_0)
  +\bar{x}_0^T\Pi \bar{x}_0+2s^T(0)\bar{x}_0\big]+m.
\end{equation*}
\end{lemma}

  \emph{Proof.} From (\ref{eq63}),
  \begin{equation}\label{eq33}
    d\mathbb{E}[{z}_i]=[(A+G)\mathbb{E}[{z}_i]+B\mathbb{E}[u_i]+f]dt, \ \mathbb{E}[{z}_i](0)=x_{i0}.
    \end{equation}
   Applying It\^{o}'s formula to $\|z_i-\mathbb{E}[{z}_i]\|^2_P$, we have
  \begin{align}\label{eq48}
    &\mathbb{E}\big[\|z_i(T)-\mathbb{E}[z_i(T)]\|^2_P-\|x_{i0}-\bar{x}_{0}\|^2_{P}\big]\cr
    =&\mathbb{E}\int_0^T\Big\{2\big\langle z_i-\mathbb{E} [z_i],P[A(z_i-\mathbb{E} [z_i])+B(u_i-\mathbb{E} [u_i])\big\rangle
    +\|Cz_i+Du_i+\sigma\|^2_P   \Big\}dt\cr
   =&\mathbb{E}\int_0^T\Big\{\big\langle(A^TP+PA+C^TPC)(z_i-\mathbb{E} [z_i]),z_i-\mathbb{E} [z_i]\big\rangle\cr
   &+2\big\langle(B^TP+D^TPC)(z_i-\mathbb{E} [z_i]), u_i-\mathbb{E} [u_i]\big\rangle \cr
   &+\!\big\langle\!u_i-\mathbb{E}[u_i],D^TPD(u_i-\mathbb{E}[u_i])\!\big\rangle\!+\!\big\langle\!\mathbb{E} [u_i],D^TPD\mathbb{E}[u_i]\!\big\rangle+\langle\sigma,P \sigma\rangle\cr
   &+\big\langle C^TPC\mathbb{E} [z_i]+2C^TP\sigma,\mathbb{E} [z_i]\big\rangle
   +2\big\langle D^TPC\mathbb{E} [z_i]+D^TP\sigma,\mathbb{E} [u_i]\big\rangle
 \Big\}dt.
  \end{align}
   From (\ref{eq14}) and (\ref{eq33}),
    \begin{align}\label{eq49}
 \mathbb{E}[z_i(T)]^T\Pi\mathbb{E}[z_i(T)]-\bar{x}_{0}^T{\Pi}\bar{x}_{0}
%=& \mathbb{E}\int_0^T\Big\{ 2\big\langle\Pi\mathbb{E}[z_i],(A+G)\mathbb{E}[z_i]+B\mathbb{E}[u_i]+f\big\rangle\cr
%&    +\big|C_0\mathbb{E}[z_i]+D_0\mathbb{E}[u_i]+\sigma_0\big|^2_{\Pi}
%+2\mathbb{E}[z_i]^TM[C_0\mathbb{E}[z_i]+D_0\mathbb{E}[u_i]+\sigma_0]\Big\}dt\cr
=&\mathbb{E}\int_0^T\big\{\langle[(A+G)^T\Pi+\Pi(A+G)]\mathbb{E}[z_i],\mathbb{E}[z_i]\rangle\cr
&+2\langle B^T\Pi\mathbb{E}[z_i],\mathbb{E}[u_i]\rangle  +2\langle \Pi f, \mathbb{E}[z_i]\rangle \big\}dt.
\end{align}
Also, applying It\^{o}'s formula to $\langle s,\mathbb{E}[z_i]\rangle$, we have
\begin{align}\label{eq50}
     \mathbb{E}[z_i(T)]^Ts(T)-\bar{x}_{0}^T s (0)
= & \mathbb{E}\int_0^T \Big\{ \langle-\big[A+G-B\Upsilon^{\dag}\big(B^T\Pi+D^TPC\big)\big]^Ts, \mathbb{E}[z_i]\rangle\nonumber\\
&-\langle\big[C-D\Upsilon^{\dag}\big(B^T\Pi+D^TPC\big)\big]^TP\sigma+\Pi f-\bar{\eta}, \mathbb{E}[z_i]\rangle\big\}\nonumber\\
&+ \langle(A+G)\mathbb{E}[z_i]+B\mathbb{E}[u_i]+f,s\rangle\Big\}dt\cr
= & \mathbb{E}\!\!\int_0^T\!\! \Big\{\big\langle (\Pi B+C^TPD)\Upsilon^{\dag}(B^Ts+D^TP\sigma),\mathbb{E}[z_i]\big\rangle +\langle s,f\rangle\nonumber\\
&-\langle C^TP\sigma+\Pi f-\bar{\eta},\mathbb{E}[z_i]\rangle+\langle B^Ts,\mathbb{E}[u_i]\rangle\Big\}dt.
     \end{align}
   Denote $\Psi\stackrel{\Delta}{=}B^TP+D^TPC$. By (\ref{eq48})-(\ref{eq50}), we obtain
  \begin{align*}
    \mathcal{J}_i(u_i)
    =&\mathbb{E}\int_0^{\infty}(\|z_i-\Gamma \mathbb{E}[{z}_i]-\eta\|^2_Q+\|u_i\|^2_R)dt\cr
    =&\mathbb{E}\int_0^{\infty}\Big[\|z_i-\mathbb{E}[{z}_i]\|^2_Q+\|(I-\Gamma)\mathbb{E}[{z}_i]\|^2_Q-2\bar{\eta}^T\mathbb{E}[{z}_i]+\|\eta\|^2_Q\cr
    &+\|u_i-\mathbb{E}[u_i]\|^2_R+\|\mathbb{E}[u_i]\|^2_R\Big]dt\cr
    =&\mathbb{E}\big[\|x_{i0}-\bar{x}_0\|^2_P+\bar{x}_0^T\Pi \bar{x}_0+2s^T(0)\bar{x}_0\big]-\lim_{T\to \infty} \mathbb{E}\big\{\|z_i(T)-\mathbb{E}[z_i(T)]\|^2_P\cr
    &+\mathbb{E}[z_i(T)]^T{\Pi}\mathbb{E}[z_i(T)]+2z_i(T)^Ts(T)\big\}+\mathbb{E}\!\!\int_0^{\infty}\!\!\Big[\langle\Psi^T\Upsilon^{\dag}\Psi(z_i-\mathbb{E}[{z}_i]), z_i-\mathbb{E}[{z}_i]\big\rangle\cr
    &+2\langle \Psi,z_i-\mathbb{E}[{z}_i]\rangle+\langle\Upsilon( u_i\!-\!\mathbb{E}[u_i]),  u_i\!-\!\mathbb{E}[u_i]\rangle\cr
    &   + \langle(B^T\Pi+D^TPC)^T\Upsilon^{\dag}(B^T\Pi+D^TPC)\mathbb{E}[{z}_i],\mathbb{E}[{z}_i]\rangle\cr
&+2\langle(B^T\Pi+D^TPC)\mathbb{E}[{z}_i]+B^Ts+D^TP\sigma,\mathbb{E}[u_i]\rangle+\langle\Upsilon \mathbb{E}[u_i],\mathbb{E}[u_i]\rangle\cr
 &+ \langle(B^T\Pi+D^TPC)^T\Upsilon^{\dag}(B^Ts+D^TP\sigma),\mathbb{E}[{z}_i]\rangle+2\langle s,f\rangle+\langle P\sigma, \sigma\rangle+\langle Q\eta,\eta\rangle\Big]dt\cr
    =&\mathbb{E}\big[\|x_{i0}-\bar{x}_0\|^2_P+\bar{x}_0^T\Pi \bar{x}_0+2s^T(0)\bar{x}_0\big]+\mathbb{E}\int_0^{\infty}\Big[\big\|u_i-\mathbb{E}[u_i]
    +\Upsilon^{\dag}\Psi(z_i-\mathbb{E}[{z}_i])\big\|^2_{\Upsilon}\cr
    &+\big\|\mathbb{E}[u_i]+\Upsilon^{\dag}(B^T\Pi+D^TPC)\mathbb{E}[{z}_i]+B^Ts+D^TP\sigma\big\|^2_{\Upsilon}\Big]dt+m\cr
    \geq&\mathbb{E}\big[\|x_{i0}-\bar{x}_0\|^2_P+\bar{x}_0^T\Pi \bar{x}_0+2s^T(0)\bar{x}_0\big]+m.
  \end{align*}
$\hfill \Box$

\begin{lemma}\label{lem5.2}
Let A1), A4)-A7) hold. Then
$$\mathbb{E}\int_0^{\infty}\|\hat{x}_i-\hat{z}_i\|^2dt=O\Big(\frac{1}{N}\Big),$$
\end{lemma}
where $\hat{z}_i$ is the closed-loop state of $z_i$ in (\ref{eq63}).

\emph{Proof.} After applying the control (\ref{eq64}) into the dynamics (\ref{eq63}), we have
\begin{equation*}
\begin{aligned}
  d\hat{z}_i=&\big[A\hat{z}_i\!-\!B\Upsilon^{\dag}[(B^TP\!+\!D^TPC)\hat{z}_i\!+\!B^T(\Pi\!-\!P)\mathbb{E}[\hat{z}_i]
  +B^Ts+D^TP\sigma]+G\mathbb{E}[\hat{z}_i]+f\big]dt\cr
  &+\!\big[C\hat{z}_i\!-\!D\Upsilon^{\dag}[(B^TP\!+\!D^TPC)\hat{z}_i\!+\!B^T(\Pi\!-\!P)\mathbb{E}[\hat{z}_i]
  +B^Ts+D^TP\sigma]+\sigma\big]dW_i,
\end{aligned}
\end{equation*}
  which leads to
\begin{equation*}
  d\mathbb{E}[\hat{z}_i]=[(A+G-(B\Upsilon^{\dag}B^T\Pi+D^TPC))\mathbb{E}[\hat{z}_i]+f]dt,\ \mathbb{E}[\hat{z}_i(0)]=\bar{x}_{0}.
\end{equation*}
By comparing this with (\ref{eq18}), we
can verify that $\mathbb{E}[\hat{z}_i]=\bar{x}.$ From (\ref{eq20i}),
\begin{equation*}
d(\hat{x}_i-\hat{z}_i)=\bar{A}(\hat{x}_i-\hat{z}_i)dt+G(\hat{x}^{(N)}-\mathbb{E}[\hat{z}_i])dt
+\bar{C}(\hat{x}_i-\hat{z}_i)dW_i.
\end{equation*}
This implies
 $$\hat{x}_i(t)-\hat{z}_i(t)=\int_0^t\Phi_i(t-\tau)G[\hat{x}^{(N)}(\tau)-\mathbb{E}[\hat{z}_i(\tau)]]d\tau,$$
where $\Phi_i$ satisfies
$$d\Phi_i(t)=\bar{A}\Phi_i(t)dt+\bar{C}\Phi_i(t)dW_i,\ \Phi_i(t)=I.$$
By Schwarz's inequality and Theorem \ref{thm4},
\begin{align*}
 \mathbb{E}\int_0^{\infty} \|\hat{x}_i(t)-\hat{z}_i(t)\|^2dt
 =&\mathbb{E}\int_0^{\infty}\Big\|\int_0^t\Phi_i(t-\tau) G(\hat{x}^{(N)}(\tau)-\mathbb{E}[\hat{z}_i(\tau)])d\tau\Big\|^2dt\\
 \leq&\mathbb{E}\int_0^{\infty}t\int_0^t\big\|\Phi_i(t-\tau)\|^2\| G(\hat{x}^{(N)}(\tau)-\mathbb{E}[\hat{z}_i(\tau)])\big\|^2d\tau dt\\
 =&\mathbb{E}\int_0^{\infty}\big\|G(\hat{x}^{(N)}(\tau)-\mathbb{E}[\hat{z}_i(\tau)])\big\|^2\int_{\tau}^{\infty}t\|\Phi_i(t-\tau) \|^2 dt d\tau\cr
 \leq & c \mathbb{E}\int_0^{\infty} \big\|\hat{x}^{(N)}(\tau)-\mathbb{E}[\hat{z}_i(\tau)]\big\|^2 d\tau=O(\frac{1}{N}).
\end{align*}
$\hfill \Box$

\emph{Proof of Theorem \ref{thm5.1}.}
Note that
$\mathbb{E}[\hat{z}_i]=\bar{x}$. We have
\begin{align*}
\frac{1}{N}J_{\rm soc}(\hat{u})=&\frac{1}{N}\sum_{i=1}^N\mathbb{E}\int_0^{\infty}\Big[\|\hat{x}_i-\Gamma \hat{x}^{(N)}+\eta)\|^2_Q\cr
&+\|\Upsilon^{\dag}[(B^TP+D^TPC)\hat{x}_i
+B^T(\Pi-P)\bar{x}+B^Ts+D^TP\sigma]\|^2_R\Big]dt\cr
=&\frac{1}{N}\sum_{i=1}^N\mathbb{E}\int_0^{\infty}\Big[\|\hat{z}_i-\Gamma (\mathbb{E}[\hat{z}_i]+\eta)
+\hat{x}_i-\hat{z}_i+\Gamma \hat{x}^{(N)}-\Gamma \mathbb{E}[\hat{z}_i]\|^2_Qdt\cr
&+\|\Upsilon^{\dag}[(B^TP+D^TPC)(\hat{z}_i+\hat{x}_i-\hat{z}_i)
+B^T(\Pi-P) \mathbb{E}[\hat{z}_i]+B^Ts+D^TP\sigma]\|^2_R\Big]dt.
\end{align*}

By Schwarz's inequality, and Lemma \ref{lem5.2}, one can obtain
\begin{align*}
&|\frac{1}{N}J_{\rm soc}(\hat{u})-\frac{1}{N}\mathcal{J}_{\rm soc}(\hat{u})|\cr
\leq&\frac{1}{N}\sum_{i=1}^N \mathbb{E}\int_0^{\infty}\big[\|\hat{x}_i-\hat{z}_i\|^2_Q+\|\Gamma( \hat{x}^{(N)}-\mathbb{E}[\hat{z}_i])\|^2_Q\big]dt\cr
&+\frac{c_1 }{N}\sum_{i=1}^N\Big( \mathbb{E}\int_0^{\infty}\|\hat{x}_i-\hat{z}_i\|^2_Qdt\Big)^{1/2}
+\frac{c_2 }{N}\sum_{i=1}^N
\Big(\mathbb{E}\int_0^{\infty}\|\Gamma( \hat{x}^{(N)}-\mathbb{E}[\hat{z}_i])\|^2_Qdt\Big)^{1/2}\cr
%=&\mathbb{E}\big[\|x_{i0}-\bar{x}_0\|_P^2+\|\bar{x}_0\|_{\Pi}^2+2s^T(0)\bar{x}_0\big]\cr
\leq& O(1/\sqrt{N}).
\end{align*}
From this and Lemma \ref{lem5.1}, the theorem follows.
$\hfill \Box$

\section{Numerical Example}
In this section, a numerical example is given to illustrate the effectiveness of the proposed decentralized control laws.

\begin{figure}[H]
	\centering
	\includegraphics[width=0.7\linewidth]{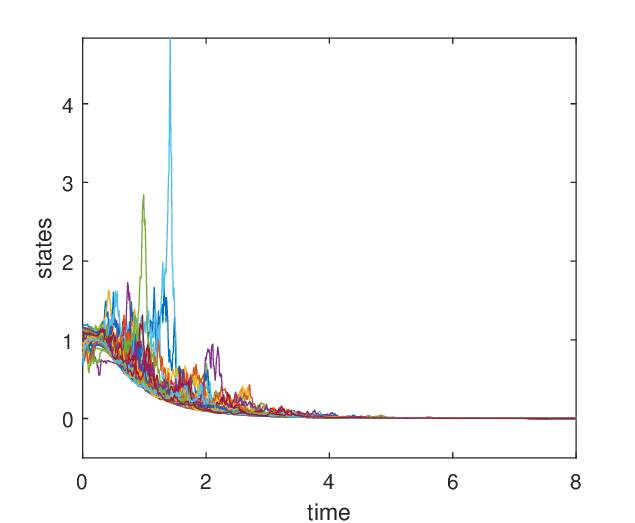}
	\caption{Curves of 30 agents.}
	\label{fig1}
\end{figure}

\begin{figure}[H]
	\centering
	\includegraphics[width=0.65\linewidth]{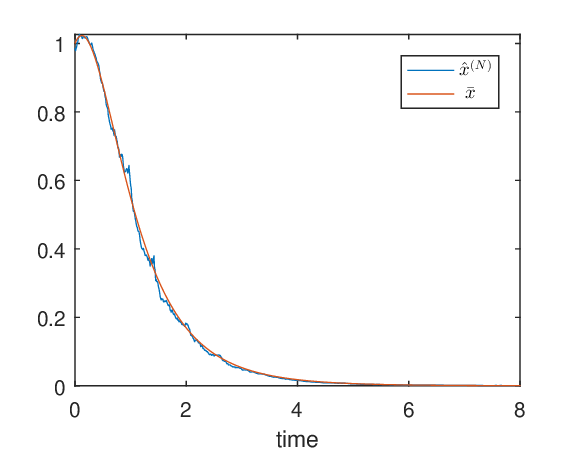}
	\caption{Curves of $\hat{x}^{(N)}$ and $\bar{x}$.}
	\label{fig2}
\end{figure}

We consider a scalar system with $50$ agents in Problem (P0). Take $A=0.1, B=C=D=Q=1,R =-0.2,G=-0.1,f=e^{-t},\eta=\frac{1}{t+1},\sigma=0.1$, and $\Gamma=-0.2 $. The initial states of $50$ agents are taken independently from a normal distribution $N(1,0.1)$. The Riccaiti equations (\ref{eq15})-(\ref{eq16}) admit solutions $P=0.6808$ and $\Pi=0.3290$, respectively.
  Then, under the control law (\ref{eq14}), the state trajectories of agents are shown in Fig. \ref{fig1}. After the transient phase, the states of agents achieve an agreement. The trajectories of $\bar{x}$ and $\hat{x}^{(N)}$ in (P0) are shown in Fig. \ref{fig2}.
It can be seen that $\bar{x}$ and $\hat{x}^{(N)}$ coincide well, which illustrates the consistency of mean field approximations.

\begin{figure}[H]
	\centering
	\includegraphics[width=0.65\linewidth]{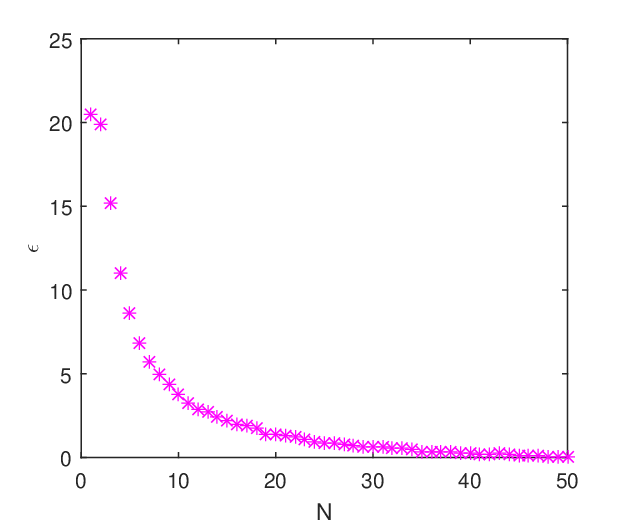}
	\caption{Curves of $\varepsilon$ with resect to $N$.}
	\label{varepsilon}
\end{figure}

The cost gap $\varepsilon$ between centralized and decentralized optimal controls is demonstrated in Fig. \ref{varepsilon} where the agent number $N$ grows from 1 to 50.

\section{Concluding Remarks}

In this paper, we have considered uniform stabilization and asymptotic optimality for indefinite mean field LQ social control systems with multiplicative noises. %(teams)
By decoupling FBSDEs, we design the decentralized control laws, which are further shown to be asymptotically optimal. Some equivalent conditions are further given for uniform stabilization of all the subsystems.

The interesting generalization is to consider mean field LQ control systems with model uncertainty by handling FBSDEs. Also, the variational analysis may be applied to leader-follower models to construct decentralized social control.

\appendices

\section{Proof of Theorem \ref{thm1}}\label{app1}
\def\theequation{A.\arabic{equation}}
\setcounter{equation}{0}
\emph{Proof of Theorem \ref{thm1}}.
(i) Suppose that $\check{u}_i$ satisfies $R\check{u}_i+B^Tp_i+D^T\beta_i^i=0,$ where $\{p_i, \beta_i^j,i,j=1,\cdots,N\}$ is a set of solutions to the equation system
\begin{align}\label{eq6}
dp_i=\alpha_idt+\beta_i^idW_i+\sum_{j\not =i}\beta_i^jdW_j,
p_i(T)=p_{iT},\ i=1,\cdots,N.
\end{align}
Here $\alpha_i, \beta_i^j, p_{iT}$, $i,j=1,\cdots,N$ are to be determined. %$\check{f}\in {\mathcal U}_{c}^{\rm F}$ is a minimizer to Problem (P1$^\prime$).
 Denote by $\check{x}_i$ the state of agent $i$ under the control $\check{u}_i$. For any $u_i\in L^2_{{\mathcal F}}(0, T; \mathbb{R}^r) $ and $\theta\in \mathbb{R}\ (\theta \not= 0)$, let $u_i^{\theta}=\check{u}_i+\theta u_i$. Denote by $x_i^{\theta}$ the solution of the following perturbed state equation:
$$ \begin{aligned}
  &dx_i^{\theta}=\big(Ax_i^{\theta}+B(\check{u}_i+\theta u_i)+\frac{G}{N}\sum_{i=1}^Nx^{\theta}_i\big)dt
  +(Cx_i^{\theta}+Du_i^{\theta}+\sigma) dW_i,\cr
 &x_i^{\theta}(0)=x_{i0},\ i=1,2,\cdots,N.
 \end{aligned}$$ Let $y_i=(x_i^{\theta}-\check{x}_i)/\theta$. %$y^{(N)}=\sum_{i=1}^Ny_i/{N}$ and
%$\textbf{y}=[y_1^T,\cdots,y_N^T]$.
It can be verified that
 $y_i$ satisfies (\ref{eq4aa}).
%By (\ref{eq1}),
%\begin{equation}\label{eq4}
%  d(\delta x_i)=[A(\delta x_i)+(\delta x^{\theta})_i+\delta f]dt,\quad \delta x_i(0)=0, \quad i=1,2,\cdots,N.
%\end{equation}
Then by It\^{o}'s formula, for any $i=1,\cdots,N$,
\begin{align*}
\mathbb{E}[\langle p_{iT},y_i(T)\rangle]
= &\mathbb{E}[\langle p_i(T),y_i(T)\rangle-\langle p_i(0),y_i(0)\rangle]\cr
= &\mathbb{E}\int_0^T \big[\langle \alpha_i,y_i\rangle+\langle p_i,Ay_i+Gy^{(N)}+Bu_i\rangle
+\langle\beta_i^i,Cy_i+Du_i\rangle\big] dt,
\end{align*}
which implies
\begin{align}\label{eq7}
&\sum_{i=1}^N\mathbb{E}[\langle p_{iT},y_i(T)\rangle]\cr
=&\sum_{i=1}^N \mathbb{E}\int_0^T   \big[\langle \alpha_i,y_i\rangle+\langle p_i,Ay_i+Gy^{(N)}+Bu_i\rangle
+\langle\beta_i^i,Cy_i+Du_i\rangle\big] dt\cr
= &\sum_{i=1}^N \mathbb{E}\int_0^T   \Big[\langle \alpha_i+A^Tp_i+C^T\beta_i^i,y_i\rangle+\langle B^Tp_i
+D^T\beta_i^i,u_i\rangle \Big]dt+\mathbb{E}\int_0^T\big\langle \sum_{i=1}^Np_i,\frac{G}{N} \sum_{i=1}^Ny_i\big\rangle dt\cr
=& \sum_{i=1}^N \mathbb{E}\int_0^T   \big[\langle \alpha_i+A^Tp_i+G^Tp^{(N)}+C^T\beta_i^i,y_i\rangle
+\langle B^Tp_i+D^T\beta_i^i,u_i\rangle\big] dt.
\end{align}
From (\ref{eq3}), we have
\begin{equation}\label{eq5a}
\begin{aligned}
 &\check{J}_{\rm soc}^{\rm F}(\check{u}+\theta u)-\check{J}_{\rm soc}^{\rm F}(\check{u})=2\theta I_1+{\theta^2}I_2
\end{aligned}
\end{equation}
where $\check{u}=(\check{u}_1,\cdots,\check{u}_N)$, and
\begin{align*}
I_1\stackrel{\Delta}{=}&\sum_{i=1}^N\mathbb{E}\Big\{\int_0^T  \big[\big\langle Q\big(\check{x}_i-(\Gamma\check{x}^{(N)}+\eta)\big),y_i-\Gamma y^{(N)}\big\rangle
+\langle R \check{u}_i,u_i\rangle \big]dt\cr
&+\big\langle H\big(\check{x}_i(T)-(\Gamma_0\check{x}^{(N)}(T)+\eta_0)\big),
y_i(T)-\Gamma_0 y^{(N)}(T)\big\rangle\Big\},\cr
I_2\stackrel{\Delta}{=}&\sum_{i=1}^N\mathbb{E}\int_0^T \big[\big\|y_i
   -\Gamma y^{(N)}\big\|^2_{Q}
+\|u_i\|^2_{R}\big]dt
+\sum_{i=1}^N\mathbb{E}\|y_i(T)-\Gamma_0 y^{(N)}(T)\|^2_H.
\end{align*}
Note that
\begin{align*}
  &\sum_{i=1}^N\mathbb{E}\Big\{\int_0^T\big\langle Q\big(\check{x}_i-(\Gamma\check{x}^{(N)}+\eta)\big),\Gamma y^{(N)}\big\rangle dt
  +\big\langle H\big(\check{x}_i(T)-(\Gamma_0\check{x}^{(N)}(T)+\eta_0)\big),\Gamma_0 y^{(N)}(T)\big\rangle\Big\}\cr
 % =&
%\mathbb{ E}\int_0^T \Big\langle \Gamma^TQ  \sum_{i=1}^N\big(\check{x}_i-(\Gamma\check{x}^{(N)}+\eta)\big),\frac{1}{N}  \sum_{j=1}^Ny_j\Big\rangle  dt\cr
 =& \sum_{j=1}^N \mathbb{E}\Big\{\int_0^T \Big\langle  \frac{\Gamma^TQ}{N} \sum_{i=1}^N\big(\check{x}_i-(\Gamma\check{x}^{(N)}+\eta)\big), y_j\Big\rangle  dt\cr
 &+\Big\langle  \frac{\Gamma^T_0H}{N} \sum_{i=1}^N\big(\check{x}_i(T)-(\Gamma_0\check{x}^{(N)}(T)+\eta_0)\big), y_j(T)\Big\rangle\Big\}\cr
 =& \sum_{j=1}^N \mathbb{E}\Big\{\int_0^T\big\langle  {\Gamma^TQ} \big((I-\Gamma)\check{x}^{(N)}-\eta\big), y_j\big\rangle  dt
 + \big\langle  {\Gamma^T_0H} \big((I-\Gamma_0)\check{x}^{(N)}(T)-\eta_0\big), y_j(T)\big\rangle\Big\}.
\end{align*}
From (\ref{eq7}), one can obtain that
\begin{align}\label{eq10d}
 I_1=&\sum_{i=1}^N\mathbb{E}\int_0^T\Big[\big\langle Q\big(\check{x}_i-(\Gamma\check{x}^{(N)}+\eta)\big),y_i-\Gamma y^{(N)}\big\rangle
 +\langle R\check{u}_i+ B^Tp_i+D^T\beta_i^i,u_i\rangle \Big]dt\cr
 &+\sum_{i=1}^N\mathbb{E}\big[\big\langle H\big(\check{x}_i(T)-(\Gamma_0\check{x}^{(N)}(T)+\eta_0)\big),
 y_i(T)-\Gamma_0 y^{(N)}(T)\big\rangle-\langle p_{iT},y_i(T)\rangle\big]\cr
 &+\sum_{i=1}^N \mathbb{E}\int_0^T \langle \alpha_i+A^Tp_i+G^Tp^{(N)}+C^T\beta_i^i,y_i\rangle dt\cr
 =&\sum_{i=1}^N\mathbb{E}\int_0^T\Big\langle R\check{u}_i+B^Tp_i+D^T\beta_i^i,u_i\Big\rangle dt\cr
 &+\sum_{i=1}^N\mathbb{E}\Big\{\int_0^T\Big\langle
 Q\check{x}_i-Q_{\Gamma}\check{x}^{(N)}-\bar{\eta}
 %&-{\Gamma^TQ} \big((I-\Gamma)\check{x}^{(N)}-\eta\big)\cr
 +\alpha_i+A^Tp_i+G^Tp^{(N)}+C^T\beta_i^i, y_i\Big\rangle dt\cr
 &+\big\langle H \check{x}_i(T)- H_{\Gamma_0} \check{x}^{(N)}(T)-\bar{\eta}_0-p_{iT}, y_i(T)\big\rangle\!\Big\}\!.
\end{align}
From (\ref{eq5a}), $\check{u}=(\check{u}_1,\cdots,\check{u}_N)$ is a minimizer to Problem (P1) if and only if
$I_2\geq0$ and $I_1=0 $.
By Proposition \ref{prop1}, $I_2\geq0$ if and only if (P1) is convex. $I_1=0 $ is equivalent to
$$\left\{
\begin{aligned}
&\alpha_i=-\big[A^Tp_i+C^T\beta_i^i-Q\check{x}_i-Q_{\Gamma}\check{x}^{(N)}-\bar{\eta}+G^Tp^{(N)}\big],\cr
&p_{iT}=H \check{x}_i(T)- H_{\Gamma_0} \check{x}^{(N)}(T)-\bar{\eta}_0,\cr
&R\check{u}_i+B^Tp_i+D^T\beta_i^i=0.
\end{aligned}
\right.
$$
Thus, we have the following optimality system:
\begin{equation}\label{eq8}
\left\{
\begin{aligned}
&d\check{x}_i=(A\check{x}_i+B\check{u}_i+G\check{x}^{(N)}+f)dt
+(C\check{x}_i+D\check{u}_i+\sigma) dW_i,\\
&d\check{p}_i=-[A^T\check{p}_i\!+\!G^T\check{p}^{(N)}\!+\!C^T\check{\beta}_i^i\!+\!Q\check{x}_i\!-\!Q_{\Gamma}\check{x}^{(N)}\!-\!\bar{\eta}]dt
+\sum_{j=1}^N\check{\beta}_i^jdW_j,\cr
&R\check{u}_i+B^T\check{p}_i+D^T\check{\beta}_i^i=0,\ i=1,\cdots,N.\\
&\check{x}_i(0)={x_{i0}}, \quad \check{p}_i(T)=H \check{x}_i(T)- H_{\Gamma_0} \check{x}^{(N)}(T)-\bar{\eta}_0.
\end{aligned}\right.
\end{equation}
This implies that FBSDE (\ref{eq4a}) admits a solution $(\check{x}_i,\check{p}_i,\check{\beta}_{i}^{j}, i,j=1,\cdots,N)$.

On other hand, if the equation system (\ref{eq4a}) admits a solution $(\check{x}_i,\check{p}_i, \check{\beta}_{i}^{j}, i,j =1,\cdots,N)$.
Let $\check{u}_i$ satisfy $R\check{u}_i+B^T\check{p}_i+D^T\check{\beta}_i^i=0$. If (P1) is convex, then by (\ref{eq5a}), $\check{u}$ is a minimizer to Problem (P1).
%%We proceed to show uniqueness. Suppose $(\acute{x}_i,\acute{p}_i, i=1,\cdots,N)$ is another solution of (\ref{eq4a}). Set $\acute{u}_i=-{R^{\dag}}B^T\acute{p}_i $. By (\ref{eq6}) and (\ref{eq10b}), we can show
%%$I_1=0$.
%%Since Problem (P1) is convex in $u$, then $\acute{u}=(\acute{u}_1,\cdots,\acute{u}_N)$ is the unique optimal minimizer, which implies the optimal pairs $(\acute{x}_i,\acute{p}_i, i=1,\cdots,N)$ coincides with $(x_i,p_i, i=1,\cdots,N)$. %On the contrary, if there exists another optimal control?
%%%$${J}_{soc}^T(u)\geq {J}_{soc}^T(\check{u}).$$

(ii) By Proposition \ref{prop2}, the fact that (P1) is uniformly convex implies (\ref{eq5}) admits a solution. This with \cite{YZ99} further gives FBSDE (\ref{eq4a}) admits a solution. Thus, (ii) follows. $\hfill \Box$
\section{Proof of Theorems \ref{thm3}}\label{app2}
\def\theequation{B.\arabic{equation}}
\setcounter{equation}{0}
To prove Theorem \ref{thm3}, we need two lemmas.

\begin{lemma}\label{lem1a}
	Let A1)-A3) hold. Under the control (\ref{eq15a}), we have
	\begin{equation}\label{eq23a}
	%\sup_{N\geq 1}
	\max_{0\leq t\leq T}\mathbb{E}\|\hat{x}_i(t)\|^2\leq c.
	\end{equation}	
\end{lemma}
\emph{Proof.} Let $\Phi_i(t)$
is the solution to the following stochastic differential equation:
\begin{equation}\label{eq23c}
  d\Phi_i(t)=\bar{A}\Phi_i(t)dt+\bar{C}\Phi_i(t)dW_i(t),\ \Phi_i(0)=I.
  \end{equation}
 From (\ref{eq20}),
we have
\begin{align*}
\hat{x}_i(t)\!=\Phi_i(t)x_{i0}+\Phi_i(t)\!\!\int_0^t\!\!\Phi_i^{-1}(\tau)(Gx^{(N)}(\tau)+g(\tau))d\tau
+\Phi_i(t)\int_0^t\Phi_i^{-1}(\tau)\sigma(\tau) dW_i(\tau),
\end{align*}
where $$g\stackrel{\Delta}{=}(\bar{C}^TD-B)\Upsilon^{\dag}B^T(K\bar{x}+s)+f-\bar{C}^T\sigma.$$
It can be verified that $\int_0^T\|g(t)\|^2dt\leq c$.
Note that $\mathbb{ E}\int_0^T tr[\Phi^T_i(t)\Phi_i(t)]dt<c$.
%\begin{align*}
 %\mathbb{ E}\int_0^T tr[\Phi^T_i(t)\Phi_i(t)]dt=&\mathbb{E}\int_0^T tr\big\{\exp[(\bar{A}^T+\bar{A}-\bar{C}^T\bar{C})t+(\bar{C}^T+\bar{C})W_i(t)]\big\}dt\cr
  %=&\int_0^Ttr\big\{ \exp[(A^T+\bar{A}+\bar{C}^T\bar{C})t]\big\}dt\leq C.\end{align*}
We have
\begin{align*}
 \mathbb{ E}\|\hat{x}_i(t)\|^2
 \leq & 3\mathbb{E}\|\Phi_i(t)x_{i0}\|^2+3\mathbb{E}\int_0^Ttr\big[\Phi_i^T(t-\tau)\sigma^T\sigma\Phi_i(t-\tau)\big] d\tau\\
& +3\mathbb{E}\int_0^ttr[\Phi_i^T(t-\tau)\Phi_i(t-\tau)]d\tau
\mathbb{E}\int_0^t\|G\hat{x}^{(N)}(\tau)+g(\tau)\|^2d\tau\\
\leq&c_0+ 6c_1\Big(c_2\mathbb{E}\int_0^T\frac{1}{N}\sum_{i=1}^N\|\hat{x}_i(\tau)\|^2d\tau+c_3\Big)\\
=&6c_1c_2\max_{1\leq i\leq N}\mathbb{E}\int_0^T\|\hat{x}_i(\tau)\|^2d\tau+c.
\end{align*}
By Gronwall's inequality,
$ \max_{1\leq i\leq N}\mathbb{ E}\|\hat{x}_i(t)\|^2\leq ce^{6c_1c_2t}.$ This implies (\ref{eq23a}).  \hfill$\Box$

\begin{lemma}\label{lem1}
	Let A1)-A3) hold. Under the control (\ref{eq15a}), we have
	\begin{equation}
	%\sup_{N\geq 1}
	\max_{0\leq t\leq T}\mathbb{E}\|\hat{x}^{(N)}(t)-\bar{x}(t)\|^2=O(\frac{1}{N}).
	\end{equation}	
\end{lemma}
\emph{Proof.} It follows by (\ref{eq20}) that
\begin{equation*}
d\hat{x}^{(N)}=\big[(\bar{A}+G)\hat{x}^{(N)}-B\Upsilon^{\dag}B^T(\Upsilon\bar{x}+s)+f\big]dt
+\frac{1}{N}\sum_{i=1}^N[\bar{C}\hat{x}_i-D\Upsilon^{\dag}B^T(K\bar{x}+s)
+\sigma] dW_i.
\end{equation*}
From this and (\ref{eq12a}), we have
\begin{equation*}
d(\hat{x}^{(N)}-\bar{x})=(\bar{A}+G)(\hat{x}^{(N)}-\bar{x})dt
+\frac{1}{N}\sum_{i=1}^N
[\bar{C}\hat{x}_i-D\Upsilon^{\dag}B^T(K\bar{x}+s)
+\sigma] dW_i,
\end{equation*}
which leads to
\begin{equation}\label{eq21}
\begin{aligned}
\hat{x}^{(N)}(t)-\bar{x}(t)=& e^{(\bar{A}+G)t}[\hat{x}^{(N)}(0)-\bar{x}(0)]\cr
&+\frac{1}{N}\sum_{i=1}^N\int_0^te^{(\bar{A}+G)(t-\tau)}[\bar{C}\hat{x}_i
-D\Upsilon^{\dag}B^T(K\bar{x}+s)
+\sigma] dW_i(\tau).
\end{aligned}
\end{equation}
By A1), one can obtain
\begin{align*}
&\mathbb{E}\big\|  \hat{x}^{(N)}(t)-\bar{x}(t)\big\|^2\cr
\leq\ &2\big\|e^{(\bar{A}+G)t}\big\|^2\Big\{\mathbb{E}\big\|\hat{x}^{(N)}(0)-\bar{x}_0\big\|^2
+\frac{1}{N}\int_0^t \big\|e^{-(\bar{A}+G)(t-\tau)}\big\|^2(c_1\mathbb{E}\|\hat{x}_i\|^2+c_2) d\tau\Big\}\cr
\leq \ &  \frac{2}{N}\big\|e^{(\bar{A}+G)t}\big\|^2\Big\{\max_{1\leq i\leq N}\mathbb{E}\|\hat{x}_{i0}\|^2
+c\int_0^t\big\|e^{-(\bar{A}+G)(t-\tau)}\big\|^2\big] d\tau\Big\},
\end{align*}
which completes the proof.  $\hfill \Box$

\emph{Proof of Theorem \ref{thm3}}. % 这个证明还需要一个高维凸的假设。。。
We first prove that for $u\in \mathcal{U}_c$,  $J_{\rm soc}^{\rm F}(u)< \infty$ implies that
$\mathbb{E}\int_0^T(\|x_i\|^2+\|u_i\|^2)dt<\infty,$ for all $i=1,\cdots,N$.
By A2), we have
$$\delta_0 \sum_{i=1}^N\mathbb{E}\int_0^T \|u_i\|^2dt-c\leq J_{\rm soc}^{\rm F}(u)< \infty,$$
which implies
$\sum_{i=1}^N\mathbb{E}\int_0^T\|u_i\|^2dt<c_1.$
%This leads to
%$$\mathbb{E}\int_0^T\|u^{(N)}\|^2dt\leq\frac{1}{N}\sum_{i=1}^N\mathbb{E}\int_0^T\|u_i\|^2dt<\infty,$$
%where $u^{(N)}=\frac{1}{N}\sum_{i=1}^Nu_i.$
By (\ref{eq1}) and Schwarz's inequality,
\begin{equation*}\label{eq38}
\mathbb{E}\|x_i(t)\|^2\leq c_1\mathbb{E}\int_0^t\|x^{(N)}(\tau)\|^2d\tau+c_2
\leq\frac{c_1}{N}\mathbb{E}\int_0^t\sum_{j=1}^N\|x_j(\tau)\|^2d\tau+c_2
\end{equation*}
which further gives that
$$\sum_{j=1}^N\mathbb{E}\|x_j(t)\|^2
\leq{c_1}\int_0^t\sum_{j=1}^N\mathbb{E}\|x_j(\tau)\|^2d\tau+Nc_2.$$
By Gronwall's inequality,
$$\sum_{j=1}^N\mathbb{E}\|x_j(t)\|^2\leq Nc_2e^{c_1t}\leq Nc_2e^{c_1T}.$$
Let $\tilde{x}_i=x_i-\hat{x}_i$,   $\tilde{u}_i=u_i-\hat{u}_i$ and $\tilde{x}^{(N)}=\frac{1}{N}\sum_{i=1}^N \tilde{x}_i$.
Note that it follows by Lemma \ref{lem1a} that
$$
\mathbb{E}\int_0^T\big(\|\hat{x}_i\|^2+\|\hat{u}_i\|^2)dt<\infty.
$$
Then we have
\begin{equation}\label{eq22b}
\mathbb{E}\int_0^T\big(\|\tilde{x}_i\|^2+\|\tilde{u}_i\|^2)dt<\infty.
\end{equation}
By (\ref{eq1}) and (\ref{eq20}),
\begin{equation}\label{eq32}
d\tilde{x}_i=(A\tilde{x}_i+{G}\tilde{x}^{(N)}+B\tilde{u}_i)dt+(C\tilde{x}_i+D\tilde{u}_i)dW_i, \  \tilde{x}_i(0)=0.
\end{equation}
From (\ref{eq3}), we have
\begin{align}
J^{\rm F}_{\rm soc}(u)
=\ &\sum_{i=1}^N(J_i^{\rm F}(\hat{u})+\tilde{J}_i^{\rm F}(\tilde{u})+2\hat{I}_i),
\end{align}
where
\begin{align*}
\tilde{J}_i^{\rm F}(\tilde{u})\stackrel{\Delta}{=}&\mathbb{E}\int_0^T \big[\|\tilde{x}_i-\Gamma \tilde{x}^{(N)}\|^2_Q+\|\tilde{u}_i\|^2_{R}\big]dt
+\mathbb{E}[\|\tilde{x}_i(T)-\Gamma_0(\tilde{x}^{(N)}(T))\|_H^2]\cr
\hat{I}_i=&\mathbb{E}\Big\{\!\!\int_0^T \!\!\big[\big(\hat{x}_i\!-\!\Gamma \hat{x}^{(N)}\!-\!\eta\big)^T\!Q\big(\tilde{x}_i\!-\!\Gamma \tilde{x}^{(N)}\big)\!+\!\hat{u}_i^TR\tilde{u}_i\big]dt\cr
& +\big[\hat{x}_i(T)-(\Gamma_0\hat{x}^{(N)}(T)+\eta_0)\big]^TH[\tilde{x}_i(T)-\Gamma_0 \tilde{x}^{(N)}(T)]\Big\}.
\end{align*}
By A2), $\tilde{J}_i^{\rm F}(\tilde{u})\geq 0$. We now prove $\frac{1}{N}\sum_{i=1}^N \hat{I}_i=O(\frac{1}{\sqrt{N}})$.
\begin{equation}\label{eq24a}
\begin{aligned}
\sum_{i=1}^N \hat{I}_i
%=\ &\sum_{i=1}^N 2\mathbb{E}\int_{0}^{T}\Big\{\tilde{x}_i^T\big[Q(\hat x_i-\Gamma \hat{x}^{(N)}-\eta)\\
%&-\Gamma^TQ((I-\Gamma)\hat{x}^{(N)}-\eta)\big]+\sum_{i=1}^N\hat{u}_i^TR\tilde{u}_i\Big\}dt\cr
=& \sum_{i=1}^N \mathbb{E}\!\int_{0}^{T}\!\! \Big\{\tilde{x}_i^T\big(Q\hat x_i-Q_{\Gamma}\bar{x}-\bar{\eta}\big)+\hat{u}_i^TR\tilde{u}_i\Big\}dt
+\sum_{i=1}^N\mathbb{E}\int_{0}^{T} (\hat{x}^{(N)}-\bar{x})^TQ_{\Gamma}\tilde{x}_idt\cr
&+\sum_{i=1}^N\mathbb{E}\big[\tilde{x}_i^T(T)(H\hat x_i(T)-H_{\Gamma_0}\bar{x}(T)-\bar{\eta}_0)
+(\hat{x}^{(N)}(T)-\bar{x}(T))^TH_{\Gamma_0}\tilde{x}_i(T)\big].
\end{aligned}
\end{equation}
%where $\xi=\hat{x}^{(N)}-\bar{x}$.
Denote $\hat{p}_i(t)=P\hat{x}_i(t)+K\bar{x}(t)+s(t)$.
Then by (\ref{eq8a})-(\ref{eq10a}) and It\^{o}'s formula,
\begin{align}\label{eq31}
  d\hat{p}_i=&-\big[A^TP+PA+C^TPC+Q-\big(B^TP+D^TPC\big)^T
  \Upsilon^{\dag}\big(B^TP+D^TPC\big)\big]\hat{x}_idt\cr
  &+P\big[\bar{A}\hat{x}_i-B\Upsilon^{\dag}(B^T
  (K\bar{x}+s)+D^TP\sigma)+G\hat{x}^{(N)}+f\big]dt\cr
  &+P[\bar{C}\hat{x}_i-D\Upsilon^{\dag}(B^T(K\bar{x}+s)+D^TP\sigma)+\sigma] dW_i\cr
  &-\big[(A+G)^TK+K(A+G)-(B^TP+D^TPC)^T
  \Upsilon^{\dag}B^TK-KB\Upsilon^{\dag}B^TK+G^TP+PG\cr
  &-KB\Upsilon^{\dag}(B^TP+D^TPC)-Q_{\Gamma}\big]\bar{x}dt+K\big\{(A+G)\bar{x}-B{\Upsilon^{\dag}}[B^T(P+K)+D^TPC]\bar{x}\cr
  &- B{\Upsilon^{\dag}}(B^Ts
  +D^TP\sigma)+f\big\}dt-\Big\{\big[A+G-B\Upsilon^{\dag}\big(B(P+K)+D^TPC\big)\big]^Ts\cr
  &+(P+K)f+\big[C-D\Upsilon_N^{\dag}\big(B(P+K)+D^TPC\big)\big]^TP\sigma-\bar{\eta}\Big\}\cr
  =&-(A^T\hat{p}_i+G^T\hat{p}^{(N)}+C^T\hat{\beta}_i^i+Q\hat{x}_i-Q_{\Gamma}\bar{x}-\bar{\eta} )dt\cr
  &+(G^TP+PG)(\hat{x}^{(N)}-\bar{x})dt+\hat{\beta}_i^idW_i,
\end{align}
where $\hat{\beta}_i^i=P(C\hat{x}_i+D\hat{u}_i+\sigma)$. By (\ref{eq15a}), we have $R\hat{u}_i=-(B\hat{p}_i+D\hat{\beta}_i^i).$
Note that $\hat{p}_i(T)=H\hat x_i(T)-H_{\Gamma_0}\bar{x}(T)-\bar{\eta}_0$.
From (\ref{eq32}) and (\ref{eq31}),
$$  \begin{aligned}
\sum_{i=1}^N\mathbb{E}\big[ \tilde{x}_i^T(T)(H\hat x_i(T)-H_{\Gamma_0}\bar{x}(T)-\bar{\eta}_0)\big]
%=\ &\sum_{i=1}^N\mathbb{E}\int_{0}^{T} \Big\{-\tilde{x}_i^T\big[ Q\hat{x}_i+G^T((P+K)\bar{x}+s)\\
%&-Q(\Gamma \bar{x}+\eta)-\Gamma^TQ\left((I-\Gamma) \bar{x}-\eta\right)\big] \\&
%+\tilde{x}_i^TPG(\hat{x}^{(N)}-\bar{x})\\
%&+ ({G}\tilde{x}^{(N)}+B\tilde{u}_i)^T(P\hat{x}_i+K\bar{x}+s)\Big\}dt\\
=&\mathbb{E}\int_{0}^{T}\sum_{i=1}^N\Big\{- \tilde{x}_i^T\big[ Q\hat{x}_i-Q_{\Gamma}\bar{x}-\bar{\eta}\big]
-\hat{u}_i^TR\tilde{u}_i\Big\}dt\cr
&+N\mathbb{E}\int_{0}^{T}(\hat{x}^{(N)}-\bar{x})^T(G^TP+PG)\tilde{x}^{(N)}dt.
\end{aligned}$$
This and (\ref{eq24a}) lead to
\begin{equation*}
\frac{1}{N}\sum_{i=1}^N \hat{I}_i=\mathbb{E}\int_{0}^{T}(\hat{x}^{(N)}\!-\bar{x})^T (Q_{\Gamma}\!+\!G^TP\!+\!PG)\tilde{x}^{(N)}dt
+\mathbb{E}[(\hat{x}^{(N)}(T)-\bar{x}(T))^TH_{\Gamma_0}\tilde{x}^{(N)}(T)\big].
\end{equation*}
By Lemma \ref{lem1}, and (\ref{eq22b}), we obtain
$$
\begin{aligned}
\Big|\frac{1}{N}\sum_{i=1}^N \hat{I}_i\Big|^2\leq &c\mathbb{E}\int_{0}^{T}\|\hat{x}^{(N)}-\bar{x}\|^2dt\cdot \mathbb{E}\int_{0}^{T}\|\tilde{x}^{(N)}\|^2dt\cr
&\times
\mathbb{E}[\|\hat{x}^{(N)}(T)-\bar{x}(T))\|^2\cdot\mathbb{E}\|\tilde{x}^{(N)}(T)\|^2,
\end{aligned}
$$
which implies $|\frac{1}{N}\sum_{i=1}^N \hat{I}_i|=O(1/\sqrt{N})$.  $\hfill \Box$
\section{Proof of Theorems \ref{thm5} and \ref{thm5b}}\label{app3}
\def\theequation{C.\arabic{equation}}
\setcounter{equation}{0}
\emph{Proof of Theorem \ref{thm5}.} (iii)$\Rightarrow$(i) was given in Theorem \ref{thm4}. We now prove (i)$\Rightarrow$(iii).
By (\ref{eq20i}),
\begin{equation}\label{eq24}
\frac{d\mathbb{E}[\hat{x}_i]}{dt}=\bar{A}\mathbb{E}[\hat{x}_i]-B\Upsilon^{\dag}B^T((\Pi-P)\bar{x}+s)
+G\mathbb{E}[\hat{x}^{(N)}]+f,\ \mathbb{E}[\hat{x}_i(0)]=\bar{x}_0.
\end{equation}
%where $\bar{A}=A-B\Upsilon^{\dag}(B^TP+D^TPC)$.
It follows from A1) that
$$\mathbb{E}[\hat{x}_i]=\mathbb{E}[\hat{x}_j]=\mathbb{E}[\hat{x}^{(N)}],\ j\not =i.$$
By comparing (\ref{eq18}) and (\ref{eq24}), we obtain
$$\frac{d(\mathbb{E}[\hat{x}_i]-\bar{x})}{dt}=(\bar{A}+G)(\mathbb{E}[\hat{x}_i]-\bar{x}),\ \mathbb{E}[\hat{x}_i(0)]-\bar{x}(0)=0,$$
which implies
\begin{equation}\label{eq24b}
  \mathbb{E}[\hat{x}_i]=\bar{x}=\mathbb{E}[\hat{x}^{(N)}].
  \end{equation}
 Note that $\|\bar{x}\|^2\leq \mathbb{E}\|\hat{x}_i\|^2$.
It follows from (\ref{eq23}) that
\begin{equation}\label{eq25}
\int_0^{\infty} \|\bar{x}(t)\|^2dt<\infty.
\end{equation}
By (\ref{eq18}), we have
\begin{equation*}
\bar{x}(t)=e^{[A+G-B\Upsilon^{\dag}(B^T\Pi+D^TPC)]t}\Big[\bar{x}_0
+\int_0^te^{-(A+G-B\Upsilon^{\dag}B^T\Pi)\tau}h(\tau)d\tau\Big] ,
\end{equation*}
where $h=-B\Upsilon^{\dag}(B^Ts+D^TP\sigma)+f$. By the arbitrariness of $\bar{x}_0$ with (\ref{eq25}) we obtain that $A+G-B\Upsilon^{\dag}(B^T\Pi+D^TPC)$ is Hurwitz. That is, $(A+G, B)$ is stabilizable.
 Note that $\mathbb{E}[x^{(N)}]^2\leq \frac{1}{N}\sum_{i=1}^N\mathbb{E}[\hat{x}_i^2]$. Then from (\ref{eq23}) we have %$q\in C_{\rho/2}([0,\infty),\mathbb{R}^n)$.
\begin{equation}\label{eq27}
\mathbb{E} \int_0^{\infty}\big\|\hat{x}^{(N)}(t)\big\|^2dt<\infty.
\end{equation}
This leads to $ \mathbb{E} \int_0^{\infty}\|k(t)\|^2dt<\infty$, where $k{=}-B\Upsilon^{\dag}[B^T((\Pi-P)\bar{x}+s)+D^TP\sigma]+G\hat{x}^{(N)}+f$.
By (\ref{eq20i}), we obtain %that $$\hat{x}_i=e^{\bar{A}t}\left(x_{i0}+\int_0^te^{-\bar{A}s}g(s)ds+\int_0^te^{-\bar{A}s}\sigma(s)dW_i(s)\right),$$
%  which implies
$$
\begin{aligned}
\mathbb{E}\|\hat{x}_i(t)\|^2=\ & \mathbb{E}\left\|\Phi_i(t)\left(x_{i0}+\int_0^t\Phi_i^{-1}(\tau)k(\tau)d\tau\right)\right\|^2,
\end{aligned}
$$
where $\Phi_i$ satisfies (\ref{eq23c}).
By (\ref{eq23}) and the arbitrariness of ${x}_{i0}$  we obtain that $\mathbb{E}\int_0^{\infty}\left\|\Phi_i(t)\right\|^2dt<\infty$, %$(\bar{A},\bar{C})$ is mean-square stable,
i.e., $[A, B; C,D]$ is stabilizable.
From (\ref{eq25}) and (\ref{eq27}),
\begin{equation}\label{eq28}
\mathbb{E} \int_0^{\infty} \big\|\hat{x}^{(N)}(t)-\bar{x}(t)\big\|^2 dt<\infty.
\end{equation}
On the other hand, it follows from (\ref{eq21d}) that
\begin{equation*}
\begin{aligned}
&\mathbb{E}\big\|\hat{x}^{(N)}(t)-\bar{x}(t)\big\|^2\\
=&\mathbb{E}\big\|e^{(\bar{A}+G)t}[\hat{x}^{(N)}(0)-\bar{x}_0]\big\|^2
+\frac{1}{N^2}\sum_{i=1}^N\mathbb{E}\int_0^t\big\|e^{(\bar{A}+G)(t-\tau)}(\bar{C}\hat{x}_i(\tau)+\bar{\sigma}(\tau)\big\|^2 d\tau.
\end{aligned}
\end{equation*}
By (\ref{eq28}) and the arbitrariness of ${x}_{i0}, i=1,\cdots,N$, we obtain that $\bar{A}+G$ is Hurwitz.

(iii)$\Rightarrow$(ii) was given in Lemma \ref{lem2a}. (ii)$\Rightarrow$(iii) was implied from \cite[Theorem 2]{LQZ19}.
\hfill{$\Box$}

\emph{Proof of Theorem \ref{thm5b}.} (iii)$\Rightarrow$(i) has been proved in Theorem \ref{thm4}.  Following (i)$\Rightarrow$(iii) of Theorem \ref{thm5}, together with \cite{AM90}, \cite{ZZC08}, we obtain
(i)$\Rightarrow$ (ii).

(ii)$\Rightarrow$(iii).  Define $V(t)=\mathbb{E}[{y}^T(t)P {y}(t)]$,
where ${y}$ satisfies (\ref{eq37a}).
%\begin{equation*}
%{d{y}}=(Ay+Bu)dt+(Cy+Du)dW(t),\quad {y}(0)={y}_0.
%\end{equation*}
Denote $V$ by $V^*$ when ${u}={u}^*(t)=-\Upsilon^{\dag}(B^TP+D^TPC)y(t)$. By (\ref{eq16}) we have
\begin{align*}
{V^*}(T)-V^*(0)
=&\mathbb{E}\Big\{y^T(t)\big[\!-Q\!-\big(B^TP\!+D^TPC\big)^T
\Upsilon^{\dag}\big(B^TP\!+D^TPC\big)\cr
&+\big(B^TP+D^TPC\big)^T
\Upsilon^{\dag}(D^TPD)
\Upsilon^{\dag}\big(B^TP+D^TPC\big) \big]y(t)\Big\}\cr
=&\mathbb{E}\big\{y^T(t)\big[-Q-\big(B^TP+D^TPC\big)^T
\Upsilon^{\dag}R\Upsilon^{\dag}\big(B^TP+D^TPC\big)\big]y(t)\big\}\\
\leq&0.
\end{align*}
Note that $V^*\geq0$. Then $\lim_{t\to\infty}V^*(t)$ exists, which implies
\begin{equation}\label{eq43}
\lim_{t_0\to\infty}[V^*(t_0)-V^*(t_0+T)]=0.
\end{equation}
Rewrite $P(t)$ in (\ref{eq8b}) by $P_{T}(t)$. Then we have $P_{T+t_0}(t_0)=P_{T}(0)$.
By (\ref{eq11}),
\begin{align*}
&\mathbb{E}\int_{t_0}^{T+t_0} [y^{T}(t)Qy(t)+{u}^T(t)R{u}(t)]dt\cr
=\ &\mathbb{E}[y^T({t_0})P_{T+t_0}(t_0)y({t_0})]\!+\! \mathbb{E}\int_0^T \!\! \big\|{u}(t)
+{\Upsilon}^{\dag}\big(B^TP_{T+t_0}(t_0)+D^TP_{T+t_0}(t_0)C\big) y(t)\big\|^2_{\Upsilon}dt\cr
\geq\ &\mathbb{E}\big\|y({t_0})\big\|^2_{P_{T+t_0}(t_0)} =\mathbb{E}\big\|y({t_0})\big\|^2_{P_{T}(0)}.
\end{align*}
This with (\ref{eq43}) implies
\begin{align*}
\lim_{t_0\to\infty}\mathbb{E}\big\|{y}({t_0})\big\|^2_{P_{T}(0)}
\leq&\lim_{t_0\to\infty}\mathbb{E} \int_{t_0}^{T+t_0} (\|{y}(t)\|_{Q}^2+\|{u}^*(t)\|^2_R)dt\\
=&\lim_{t_0\to\infty}[V^*(t_0)-V^*(t_0+T)]=0.
\end{align*}
By A5$^{\prime}$), one can obtain that there exists $T>0$ such that $P_{T}(0)>0$ (See e.g. \cite{ZQ16}, \cite{ZZC08}).
Thus, we have $\lim_{t\to\infty}\mathbb{E}\big\|\bar{y}({t})\big\|^2=0$, which implies $[A, B;C,D]$ is stabilizable.

To show that $(A+G, B)$ is stabilizable, we consider to optimize
\begin{align*}
\bar{J}(u)=&\int_0^T[\bar{y}^T(s)(C^TPC+Q-Q_{\Gamma})\bar{y}(s)
+2\bar{y}^T(s)C^TPD\bar{u}(s)+\bar{u}^T(s)\Upsilon\bar{u}(s)]ds ,
\end{align*}
where $\bar{y}$ evolve by
\begin{align}\label{eq12e}
d\bar{y}(t)=&\big[(A+G)\bar{y}(t)+B\bar{u}(t)\big]dt,\
\bar{y}(0)=\bar{y}_0.
\end{align}
Let $\bar{u}^*(t)=\Upsilon^{\dag}B^T\Pi(t)\bar{y}(t)$, where $\Pi_T(t)$ satisfies
%Denote by $\Pi_T(t)$ the solution to the equation
\begin{align}
 &\dot{\Pi}+(A+G)^T\Pi+\Pi (A+G)-\big(B^T\Pi+D^TPC\big)^T\Upsilon^{\dag}
 \big(B^T\Pi+D^TPC\big)\cr
 &+C^TPC+Q-Q_{\Gamma}=0,\ \Pi(T)=0.
\end{align}
By direct calculations,
\begin{align}\label{eq56}
  \bar{y}_0^T\Pi_T(0)\bar{y}_0
  =&\int_0^T[\bar{y}^T(s)(C^TPC+Q-Q_{\Gamma})\bar{y}(s)
  +2\bar{y}^T(s)C^TPD\bar{u}^*(s)
  +(\bar{u}^*)^T(s)\Upsilon\bar{u}^*(s)]ds\cr
  =&\int_0^T\Big[\bar{y}^T(s) \ (\bar{u}^*)^T(s)\Big]
  \left[\!\begin{array}{cc}
  Q-Q_{\Gamma}+ C^TPC&C^TPD\cr
    D^TPC& \Upsilon
    \end{array}
   \right]\!\left[\!\begin{array}{c}
     \bar{y}(s)\\
     \bar{u}^*(s)
   \end{array}\right]\!ds.
\end{align}
Note that
$$
  \left[\begin{array}{cc}
  P&PD\cr
  D^TP&R\!+\!D^TPD
  \end{array}\right]\!=\!
  \left[\begin{array}{cc}
  I&0\cr
  D^T&I
  \end{array}\right]\!\left[\begin{array}{cc}
  P&0\cr
  0&R
  \end{array}\right]\!\left[\begin{array}{cc}
  I&D\cr
  0&I
  \end{array}\right]\!.
$$
Thus, we have
$$\left[\begin{array}{cc}
P&PD\cr
D^TP&R+D^TPD
\end{array}\right]>0.$$
By Schur's lemma \cite{RCMZ01}, $P-PD\Upsilon^{\dag}D^TP\geq0$. This gives $C^T(P-PD\Upsilon^{\dag}D^TP)C\geq 0$. Using Schur's lemma again, we obtain
$$\left[\begin{array}{cc}
   C^TPC&C^TPD\cr
    D^TPC& \Upsilon
    \end{array}
   \right]\geq0. $$ %is semi-positive definite.
  Assume $\bar{y}_0^T\Pi_T(0)\bar{y}_0=0$. Then from (\ref{eq56}), we have $\int_0^T\bar{y}^T(s)(Q-Q_{\Gamma})\bar{y}(s)dt=0$,
   which implies $(I-\Gamma)\sqrt{Q}\bar{y}(s)=0, \ 0\leq s\leq T$. This together with {A5$^{\prime}$)} gives $\bar{y}_0=0$. Hence, we obtain
$\Pi_T(0)>0$.
By a similar argument as the above proof, we can obtain the stabilizability of $(A+G, B)$.
 $\hfill \Box$

\bibliographystyle{plain}

\begin{thebibliography}{10}
	
     \bibitem{AFIJ03}
	H. Abou-Kandil, G. Freiling, V. Ionescu, and G. Jank, \textit{Matrix Riccati Equations in Control and Systems Theory}. Birkhiiuser Verlag, 2003.

	\bibitem{AM90}
	B. D. O. Anderson and J. B. Moore, \textit{Optimal Control: Linear Quadratic Methods}. Englewood Cliffs, NJ: Prentice Hall, 1990.
	
	\bibitem{AM15}
	J. Arabneydi and A. Mahajan, ``Team-optimal solution of finite number of mean-field coupled LQG subsystems," in {\it Proc. 54th IEEE CDC}, Osaka, Japan, 2015, pp. 5308-5313.
	
	%\bibitem{BO82}
    %T. Basar and G. J. Olsder, \emph{Dynamic Noncooperative Game Theory}. Academic Press, London,
    %1982.
	
	%\bibitem{Afij03}
	%H. Abou-Kandil,
	%G. Freiling, V. Ionescu, And
	%G. Jank, Matrix Riccati Equations
	%In Control And Systems Theory, Birkhiiuser Verlag, 2003.
	
	\bibitem{BTN16}
	D. Bauso, H. Tembine, and T. Basar, ``Opinion dynamics in social networks through mean-field games,"  \emph{SIAM J. Control Optim.}, vol. 54, no. 6, pp. 3225-3257, 2016.
	
	%\bibitem{BSY13}
	%	A. Bensoussan, K.C.J., Sung,  and S.C.P. Yam. Linear-quadratic time-inconsistent mean field
	%	games. {\it Dynamic Games Appl.}, vol. 3, no.4, pp. 537-552, 2013.
	
	\bibitem{BSYY16}
	A. Bensoussan, K.C. Sung, S.C. Yam, and S. P. Yung, ``Linear-quadratic mean field games," \emph{J. Optimization Theory \& Applications}, vol. 169, no. 2, pp. 496-529, 2016.
	
	\bibitem{BFY13}
	A. Bensoussan, J. Frehse, and P. Yam,
	{\it Mean Field Games and Mean Field Type Control Theory}.
	Springer, New York,  2013.
	
	\bibitem{C14}
	P. E. Caines, M. Huang, and R. P. Malhame, \emph{Mean field games}, in \emph{Handbook of Dynamic Game Theory}, T. Basar and G. Zaccour Eds., Springer, Berlin, 2017.
	
	\bibitem{CS14}
	P. Chan and R. Sircar, ``Bertrand and Cournot mean field games," \emph{Applied Mathematics \& Optimization}, vol. 71, no. 3, pp. 533-569, 2015.

	\bibitem{CBM15}
	Y. Chen,  A. Busic, and S. Meyn, ``State estimation and mean field control with application to demand
	dispatch," in \emph{Proc.
		54th IEEE CDC}, Osaka, 2015, pp. 6548-6555.
	
	\bibitem{CD13}
	R. Carmona and F. Delarue, ``Probabilistic analysis of mean-field games," {\it  SIAM J. Control Optim.}, vol. 51, no. 4, pp. 2705-2734, 2013.
	
	\bibitem{CD18}
R. Carmona and F. Delarue, \textit{Probabilistic theory of mean field games with applications: I and II}. Springer-Verlag, 2018.
	
	%\bibitem{bdlp} R. Buckdahn, B. Djehiche, J. Li and S. Peng (2009). Mean-field %backward
	%stochastic differential equations: a limit approach. \emph{Annals of Probability}, %37, 1524-1565.
	
	%\bibitem{CS14} P. Chan and  R. Sircar.
	%Bertrand and Cournot mean field games, {\it Appl. Math Optim.}, vol. 71, no. 3, pp.
	%533-569, 2015.

	%\bibitem{CLZ98}
	%S. Chen, X. Li, and X. Y. Zhou.
	%Stochastic linear quadratic regulators with indefinite control weight costs. {\it SIAM J. Control Optim.}, vol. 36, no. 5, pp. 1685-1702, Sept. 1998.
	
	%\bibitem{DH15}
	%	B. Djehiche, and  M. Huang.
	%	A characterization of sub-game perfect equilibria for SDEs of
	%	mean field type. {\it Dynamic Games Appl.}, vol. 6, no. 1, pp. 55-81, 2016.
	
	%\bibitem{FJ96}
	%G. Freiling, G. Jank, Existence and comparison theorems for algebraic and continuous-time Riccati
	%differential and difference equations, \emph{J. Dynamical Control Systems}, 2, pp. 529-547, 1996.
	
	\bibitem{GS13}
	D. A. Gomes and J. Saude,
	``Mean field games models--a brief survey," {\it Dyn. Games Appl.}, vol. 4, no. 2, pp. 110-154,  2014.
	
	\bibitem{GLL11}
	O. Gu\'{e}ant, J. M. Lasry, and  P. L. Lions, ``Mean field games and applications," in \textit{Paris-Princeton Lectures on Mathematical Finance}, pp. 205-266, Springer-Verlag: Heidelberg, Germany, 2011.
	
	\bibitem{H80}
	Y. C. Ho, ``Team decision theory and information structures," in \emph{Proc. IEEE}, vol. 68, no. 6, 1980, pp. 644-654.
	
	%\bibitem{HF13}
	%M. Hu and M. Fikushima.
	%Existence, uniqueness, and computation of robust Nash equilibria in a class of multi-leader-follower games. {\it  SIAM J. Optim.}, vol. 23, no. 2,   pp. 894-916, 2013.

	\bibitem{HN16a}	
M. Huang, and S. Nguyen, ``Mean field games for stochastic
growth with relative utility," \emph{Applied Mathematics \& Optimization},
vol. 74, pp. 643-668, 2016.	

	%\bibitem{HH13}
%	J. Huang and M. Huang, ``Mean field LQG games with model uncertainty," in {\it Proc. 52nd IEEE CDC}, Florence, Italy, 2013, pp. 3103-3108.
	
	\bibitem{HH16}
	J. Huang and M. Huang, ``Robust mean field linear-quadratic-Gaussian games with model uncertainty," \emph{SIAM J. Control Optim.}, vol. 55, no. 5, pp. 2811-2840, 2017.
	
	\bibitem{H10} M. Huang, ``Large-population LQG games involving a major player: the Nash certainty equivalence principle," \emph{SIAM J. Control Optim.}, vol. 48, no.5, pp. 3318-3353, 2010.
	
	%\bibitem{hcm1} M. Huang, P. E. Caines and R. P. Malham\'{e} (2004). Uplink %power adjustment in wireless communication systems: a stochastic control% %analysis. \emph{IEEE Transactions on Automatic Control}, 49, 1693-1708.
	
    %\bibitem{HCM03}
    %M. Huang, P. E. Caines, and R. P. Malham\'e, ``Individual and mass
    %behaviour in large population stochastic wireless power control
    %problems: Centralized and Nash equilibrium solutions," in {\it Proc. 42nd IEEE CDC}, Maui, HI, 2003, pp. 98-103.
	
	\bibitem{HCM07} M. Huang, P. E. Caines, and R. P. Malham\'{e}, ``Large-population cost-coupled LQG problems with non-uniform agents: Individual-mass behavior and decentralized $\varepsilon$-Nash
	equilibria," \emph{IEEE Trans.
		Autom. Control}, vol. 52, no.9, pp. 1560-1571, 2007.
	
	\bibitem{HCM12}
	M. Huang, P. Caines, and R. Malhame, ``Social optima in mean field LQG control:
	Centralized and decentralized strategies," \emph{IEEE Trans.
		Autom. Control},
	vol. 57, no. 7, pp. 1736-1751, 2012.

	\bibitem{HMC06} M. Huang, R. P. Malham\'{e}, and P. E. Caines, ``Large population stochastic dynamic games: Closed-loop McKean-Vlasov systems and the Nash certainty equivalence principle," \emph{Communication in Information and Systems}, vol. 6, pp. 221-251, 2006.
	
	\bibitem{HN16}
	M. Huang and L. Nguyen, ``Linear-quadratic mean field teams with a major agent," in	{\it Proc. 55th IEEE CDC}, Las Vegas, NV, 2016, pp. 6958-6963,.

	%\bibitem{KOH11}
	%E. Kardes, F. Ordonez, and R.W. Hall.
	%Discounted robust stochastic games and an application
	%to queueing control. {\it Operations Research}, vol. 59, no. 2,
	%pp. 365-382, 2011.

	%\bibitem{KC13}
	%A. C. Kizilkale and P. E. Caines. Mean field stochastic adaptive control. {\it IEEE Trans. Autom. Control}, vol. 58,  no. 4,   pp. 905-920, April 2013.
	
	%\bibitem{KZ05}
	%A. J. Kurdila and M. Zabarankin. {\it Convex Functional Analysis},
	%Berlin: Birkh\"auser, 2005.
	
	%\bibitem{L84}
	%V. E.  Lambson. Self-enforcing collusion in large dynamic markets,
	%\emph{J. Econom. Theory, 34}, 282-291.
	
	\bibitem{LL07} J. M. Lasry and P. L. Lions, ``Mean field games," \emph{Japan J. Math.}, vol. 2, no. 1, pp. 229-260, 2007.

	%\bibitem{L79}
    %A. J. Laub, ``A Schur method for solving algebraic Riccati equations," \emph{IEEE Trans. Autom. Control}, vol. 24, no. 6, pp. 913-921,
    %1979.

   \bibitem{LQZ19}
    Li, H., Qi, Q., \& Zhang, H. (2019).
    Stabilization control for It\^{o} stochastic system with indefinite state and control weight costs,
    provisionally accepted by \emph{Automatica}. See https://arxiv.org/pdf/1908.07684.pdf.

	\bibitem{LZ08}
	T. Li and J.-F. Zhang, ``Asymptotically optimal decentralized control
	for large population stochastic multiagent systems," {\it IEEE Trans. Autom. Control}, vol. 53, no. 7, pp. 1643-1660, 2008.
	

	\bibitem{LW19}
Y. Liang and B.-C. Wang, ``Robust mean field social optimal control with application to opinion dynamics,"  in \emph{Proc. the 15th IEEE ICCA}, Edinburgh, 2019.


	\bibitem{LZ99}
	A. Lim and X. Y. Zhou. ``Stochastic optimal LQR control with integral quadratic constraints and indefinite control weights," {\it IEEE Trans. Autom. Control}, vol. 44, no. 7, pp. 1359-1369, 1999.

	%\bibitem{LZZ16}
	%S. Li, W. Zhang, and L. Zhao
	%On Social Optima of Non-Cooperative Mean Field Games, preprint, 2016
	
    %\bibitem{MY99}  J. Ma and J. Yong, {\it Forward-backward Stochastic Differential Equations and their
    %Applications}, Springer-Verlag, New York, 1999.
	
	\bibitem{MCH13}
	Z. Ma, D. Callaway, and I. Hiskens, ``Decentralized charging control for large populations of plug-in
	electric vehicles," \emph{IEEE Trans. Control Systems Technology}, vol. 21, no. 1, pp. 67-78, 2013.
	
    %\bibitem{M77}
    %B. P. Molinari, ``The time-invariant linear-quadratic optimal control problem," \emph{Automatica}, vol. 13, no. 4, pp. 347-357, 1977.
	
	\bibitem{MB17}
	J. Moon and T. Basar, ``Linear quadratic risk-sensitive and robust mean field games," \emph{IEEE Trans.
		Autom. Control}, vol. 62, no. 3, pp. 1062-1077, 2017.

	\bibitem{NLZ16}
Y.-H. Ni, X. Li, J.-F. Zhang, ``Indefinite mean-field stochastic linear-quadratic optimal control: From finite horizon to infinite horizon," \emph{IEEE Trans. Autom. Control}, vol. 61, no. 11, pp. 3269-3284, 2016.
	
	%\bibitem{NCMH12}
	%M. Nourian,  P. E. Caines, R. P. Malham\'e, and M. Huang.
	%Mean field control in leader-follower stochastic
	%multi-agent systems: likelihood ratio based adaptation. {\it IEEE Trans. Autom. Control},
	%vol. 57, no. 11, pp. 2801-2816, Nov. 2012.
	
    \bibitem{R62}
    R. Radner, ``Team decision problems," \emph{Annals of Mathematical Statistics}, vol. 33, no. 3,	pp. 857-881, 1962.
	
    \bibitem{QZW19}
    Q. Qi, H. Zhang, and Z. Wu, ``Stabilization control for linear continuous-time
    mean-field systems," {\it  IEEE Trans. Autom. Control}, vol. 64, no. 8, pp. 3461 - 3468, 2019.

    \bibitem{RCMZ01}
    M. A. Rami, X. Chen, J. B. Moore, and X. Y. Zhou,  ``Solvability and asymptotic behavior of generalized
    Riccati equations arising in indefinite stochastic LQ controls," \emph{IEEE Trans.
    Autom. Control}, vol. 46, no. 3, 2001.

    \bibitem{SNM18}
   R. Salhab, J. L. Ny, and R. P. Malhame, ``Dynamic collective choice: Social optima,"
\emph{IEEE Trans. Autom. Control}, vol. 63, no. 10, pp. 3487-3494, 2018.	

	%\bibitem{SY14}
	%J. Sun and J. Yong. Linear quadratic stochastic differential games: Open-loop and closed-loop saddle points, \emph{SIAM J. Control Optim.}, vol. 52, no. 6, 4082-4121, 2014.
	
    \bibitem{SLY16}
    J. Sun, X. Li, and J. Yong, ``Open-loop and closed-loop solvabilities for stochastic linear quadratic optimal control problems," \emph{SIAM J. Control Optim.}, vol. 54, no. 5, pp. 2274-2308, 2016.
	
    \bibitem{SY17}
	J. Sun and J. Yong. ``Stochastic linear quadratic optimal control problems in infinite horizon,"
	\emph{Applied Mathematics \& Optimization},
	vol. 78, pp. 145-183, 2018.
	
	\bibitem{W18}
	B.-C. Wang, ``A complete solution to mean field linear quadratic control", in \emph{Proc. 37th CCC}, Wuhan, China, 2018, pp. 1556-1563.
	
	\bibitem{WH15}
	B.-C. Wang and M. Huang, ``Mean field production output control with sticky prices: Nash and social solutions," \emph{Automatica}, vol. 100, no. 100, pp. 90-98, 2019.
	
    \bibitem{WNZ19}
    B.-C. Wang, Y.-H. Ni and H. Zhang, ``Mean field games for multi-agent systems with multiplicative noises,"  \emph{International Journal of Robust and Nonlinear Control}, in press, 2019.


	
	\bibitem{WZ13}
	B.-C. Wang and J.-F. Zhang,
	``Mean field games for large-population multiagent systems with Markov jump
	parameters," {\it SIAM J. Control Optim.}, vol. 50, no. 4,
	pp. 2308-2334, 2012.
	
	\bibitem{WZ12}
	B.-C. Wang and J.-F. Zhang, ``Distributed control of multi-agent systems with random parameters and a major agent," \emph{Automatica}, vol. 48, no. 9, pp. 2093-2106, 2012.
	
	\bibitem{WZ14}
	B.-C. Wang and J.-F. Zhang, ``Hierarchical mean field games for multiagent systems with tracking-type costs: Distributed $\varepsilon$-Stackelberg equilibria,"
	{\it  IEEE Trans. Autom. Control}, vol. 59, no. 8, pp. 2241-2247, 2014.
	
	\bibitem{WZ17}
	B.-C. Wang and J.-F. Zhang, ``Social optima in mean field linear-quadratic-Gaussian models with Markov jump parameters," \emph{SIAM J. Control Optim.}, vol. 55, no. 1, pp. 429-456, 2017.
	
	\bibitem{WZ19}
B.-C. Wang and H. Zhang, ``Indefinite linear quadratic mean field social control with multiplicative noise,"  in \emph{Proc. 15th IEEE ICCA}, Edinburgh, 2019.

	\bibitem{weintraub2008markov}
	G. Weintraub, C. Benkard, and B. Van Roy, ``Markov perfect industry dynamics with many firms,"
	\emph{Econometrica}, vol. 76, no. 6, pp. 1375--1411, 2008.
	
	%\bibitem{W68}
    %W. Wonham, ``On a matrix Riccati equation of stochastic control," \emph{SIAM J. Control Optim.}, vol. 6, no. 4, pp. 681-697, 1968.

    %\bibitem{YMMS12}
    %H. Yin, P. G. Mehta, S. P. Meyn, and U. V. Shanbhag, ``Synchronization of coupled oscillators is a game," {\it IEEE Trans. Autom. Control}, vol. 57, no. 4, pp. 920-935, April  2012.

    \bibitem{Y13}
    J. Yong, ``Linear-quadratic optimal control problems for mean-field stochastic differential equations," {\it SIAM J. Control Optim.}, vol. 51, no. 4, pp. 2809-2838, 2013.
	
	\bibitem{YZ99} J. Yong and X. Y. Zhou, \emph{Stochastic Controls: Hamiltonian Systems and HJB Equations}. Springer-Verlag, New York, 1999.
	
	\bibitem{ZX17}
	H. Zhang and J. Xu, ``Control for It\^{o} stochastic systems with input delay," \emph{IEEE Trans. Autom. Control}, vol. 62, no. 1, pp. 350-365, 2017.

	\bibitem{ZQ16}
	H. Zhang, Q. Qi, and M. Fu,
	``Optimal stabilization control for discrete-time
mean-field stochastic systems," \emph{IEEE Trans. Autom. Control}, in press, 2018.
	
	\bibitem{ZZC08}
	W. Zhang, H. Zhang, and B. S. Chen, ``Generalized Lyapunov equation approach
	to state-dependent stochastic stabilization/detectability criterion,"
	\emph{IEEE Trans. Autom. Control}, vol. 53, no. 7, pp. 1630-1642, 2008.
	
	\bibitem{ZL00}
X. Y. Zhou, and D. Li. ``Continuous-time mean-variance portfolio selection: A stochastic LQ framework." \emph{Applied Mathematics and Optimization}, vol. 42, no. 1, pp. 19-33, 2000.
\end{thebibliography}

\begin{IEEEbiography}[{\includegraphics[width=1in,height=1.25in, clip, keepaspectratio]{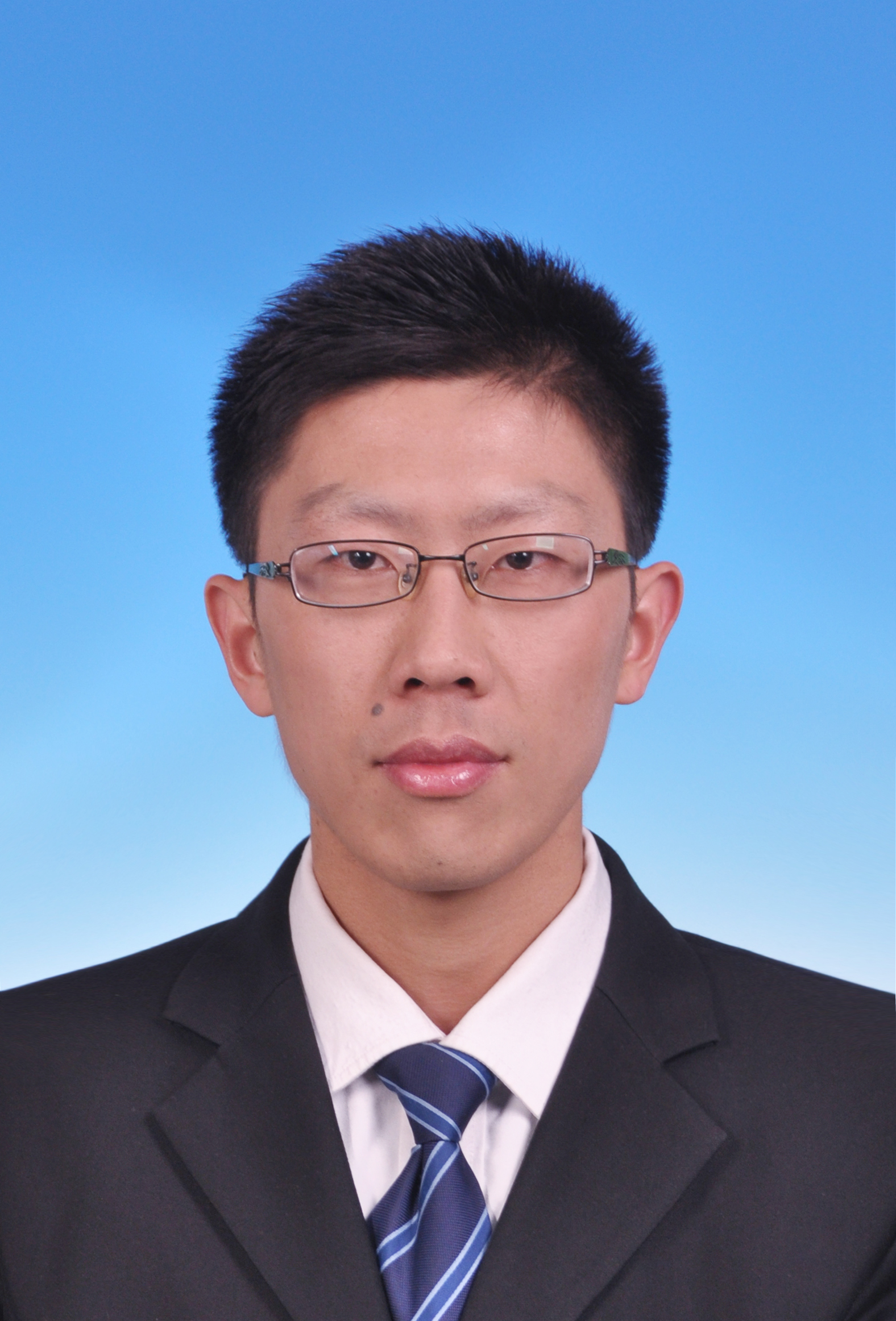}}]
{Bing-Chang Wang}
 received the M.Sc. degree in Mathematics from
Central South University, Changsha, China, in 2008, and
the Ph.D. degree in System Theory from Academy of
Mathematics and Systems Science, Chinese Academy of
Sciences, Beijing, China, in 2011. From September 2011
to August 2012, he was with Department of Electrical
and Computer Engineering, University of Alberta, Canada,
as a Postdoctoral Fellow. From September 2012 to September 2013, he was with School of Electrical Engineering and Computer
Science, University of Newcastle, Australia, as a Research
Academic.

From October 2013, he has
been with School of Control Science and Engineering, Shandong University, China, as an associate Professor. He held visiting
appointments as a Research Associate with Carleton University, Canada, from November 2014 to May 2015, and with the Hong Kong Polytechnic University from November 2016 to January
 2017. He also visited
 the Hong Kong Polytechnic University as a Research Fellow in March 2017 and May 2018.
His current research interests include mean field games, stochastic control, multiagent
systems and event based control. He received the IEEE CSS Beijing Chapter Young Author Prize in 2018.
\end{IEEEbiography}

\begin{IEEEbiography}[{\includegraphics[width=1in,height=1.25in, clip, keepaspectratio]{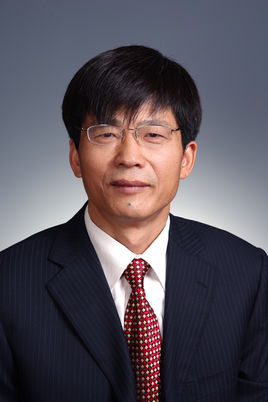}}]{Huanshui Zhang}
(SM'06) received the B.S. degree
in mathematics from Qufu Normal University, Shandong,
China, in 1986, the M.Sc. degree in control
theory from Heilongjiang University, Harbin, China,
in 1991, and the Ph.D. degree in control theory from
Northeastern University, China, in 1997.

He was a Postdoctoral Fellow at Nanyang Technological
University, Singapore, from 1998 to 2001
and Research Fellow at Hong Kong Polytechnic
University, Hong Kong, China, from 2001 to 2003.
He is currently holds a Professorship at Shandong
University, Shandong, China. He was a Professor with the Harbin Institute
of Technology, Harbin, China, from 2003 to 2006. He also held visiting
appointments as a Research Scientist and Fellow with Nanyang Technological
University, Curtin University of Technology, and Hong Kong City University
from 2003 to 2006. His interests include optimal estimation and control,
time-delay systems, stochastic systems, signal processing and wireless sensor
networked systems.
\end{IEEEbiography}

\end{document}